\documentstyle{amlts}
\begin{document}
\annalsline{155}{2002}
\received{July 10, 2000}
\startingpage{157}
\def\bye{\end{document}}
 \font\tenrm=cmr10
\input amssym.def
\input amssym.tex
 
\catcode`\@=11
\font\twelvemsb=msbm10 scaled 1100
\font\tenmsb=msbm10
\font\ninemsb=msbm10 scaled 800
\newfam\msbfam
\textfont\msbfam=\twelvemsb  \scriptfont\msbfam=\ninemsb
  \scriptscriptfont\msbfam=\ninemsb
\def\msb@{\hexnumber@\msbfam}
\def\Bbb{\relax\ifmmode\let\next\Bbb@\else
 \def\next{\errmessage{Use \string\Bbb\space only in math
mode}}\fi\next}
\def\Bbb@#1{{\Bbb@@{#1}}}
\def\Bbb@@#1{\fam\msbfam#1}
\catcode`\@=12

 \catcode`\@=11
\font\twelveeuf=eufm10 scaled 1100
\font\teneuf=eufm10
\font\nineeuf=eufm7 scaled 1100
\newfam\euffam
\textfont\euffam=\twelveeuf  \scriptfont\euffam=\teneuf
  \scriptscriptfont\euffam=\nineeuf
\def\euf@{\hexnumber@\euffam}
\def\frak{\relax\ifmmode\let\next\frak@\else
 \def\next{\errmessage{Use \string\frak\space only in math
mode}}\fi\next}
\def\frak@#1{{\frak@@{#1}}}
\def\frak@@#1{\fam\euffam#1}
\catcode`\@=12

 \newcommand{\etal}{{\it et.~al.}}
\newcommand{\ie} {{\it i.e.,\ }}
\newcommand{\eg} {{\it e.g.,\ }}
\newcommand{\cf}{{\it cf.,\ }}

\newcommand{\eqdef}{\mathbin{\stackrel{\rm def}{=}}} 
\newcommand{\R}{{\Bbb R}} 
\newcommand{\N}{{{\Bbb N}}} 
\newcommand{\Z}{{\Bbb Z}} 
\newcommand{\F}{{\Bbb F}}
\newcommand{\C}{{\Bbb C}}
\newcommand{\poly}{{{\rm poly}}} 
\newcommand{\polylog}{{{\rm polylog}}} 
\newcommand{\polyloglog}{{{\rm polyloglog}}} 
\newcommand{\loglog}{{{\rm loglog}}} 
\newcommand{\polylogloglog}{{{\rm polylogloglog}}}  
\newcommand{\zo}{\{0,1\}}
\newcommand{\suchthat}{{\;\; : \;\;}}
\newcommand{\pr}[1]{\Pr\left[#1\right]}
\newcommand{\deffont}{\em}
\newcommand{\getsr}{\gets}
\newcommand{\E}{\mathop{\mathrm E}\displaylimits}
\newcommand{\Var}{\mathop{{\rm Var}}\displaylimits}

\newcommand{\vect}[1]{{#1}}
\newcommand{\Ones}[1]{\vect{1_{#1}}}
\newcommand{\Uni}[1]{\vect{u_{#1}}}
\newcommand{\NOnes}{{\Ones{N}}}
\newcommand{\NUni}{{\Uni{N}}}

\newcommand{\PI}{{\vect{\pi}}}
\newcommand{\PIperp}{{\vect{\pi^{\perp}}}}
\newcommand{\ALPHA}{{\vect{\alpha}}}
\newcommand{\ALPHApar}{{\vect{\alpha^{\parallel}}}}
\newcommand{\ALPHAperp}{{\vect{\alpha^{\perp}}}}
\newcommand{\BETA}{{\vect{\beta}}}
\newcommand{\Clouds}{C}
\newcommand{\Unit}[1]{\vect{e}_{#1}}
\newcommand{\liftA}{{\tilde{A}}}
\newcommand{\liftB}{{\tilde{B}}}
\newcommand{\liftC}{{\tilde{C}}}
\newcommand{\Mmod}{{M'}}

\newcommand{\lamtwo}{{\overline{\lambda_2}}}

\newcommand{\Rot}[1]{{{\rm Rot}_{#1}}}
\newcommand{\GraphFam}{{\mathcal G}}
\newlength{\circlen}

\newcommand{\zigzag}{\mathbin{\raisebox{.2ex}{
      \hspace{-.4em}$\bigcirc$\hspace{-.75em}{\rm z}\hspace{.15em}}}}

\newcommand{\zigzagsub}{\mathbin{\raisebox{.1ex}{
     \scriptsize \hspace{-.5em}$\bigcirc$\hspace{-.75em}{{\ninerm z}}\hspace{.3em}}}}

\newcommand{\replacement}{\mathbin{\raisebox{.2ex}{
      \hspace{-.4em}$\bigcirc$\hspace{-.68em}{\rm r}\hspace{.15em}}}}

\newcommand{\replacementsub}{\mathbin{\raisebox{.1ex}{
     \scriptsize \hspace{-.5em}$\bigcirc$\hspace{-.70em}{{\ninerm r}}\hspace{.3em}}}}

\newcommand{\balanced}{\mathbin{\raisebox{.2ex}{
      \hspace{-.4em}$\bigcirc$\hspace{-.75em}{\lower1pt\hbox{\rm b}}\hspace{.15em}}}}

\newcommand{\balancedsub}{\mathbin{\raisebox{.1ex}{
     \scriptsize \hspace{-.5em}$\bigcirc$\hspace{-.72em}{\lower.75pt\hbox{\footnotesize b}}\hspace{.3em}}}}

\newcommand{\zigzagmod}{\mathbin{\raisebox{.2ex}{
      \hspace{-.4em}$\bigcirc'$\hspace{-1em}{\rm z}\hspace{.4em}}}}

\newcommand{\zigzagmodsub}{\mathbin{\raisebox{.1ex}{
     \scriptsize \hspace{-.5em}$\bigcirc'$\hspace{-1.1em}{{\ninerm z}}\hspace{.6em}}}}

\def\sni#1{\smallbreak\noindent{#1}. }
\def\ssni#1{\vglue-1pt\noindent\hskip18pt {#1}.}

\newcommand{\la}{{\lambda}}
\newcommand{\ra}{\to}

\newcommand{\remove}[1]{}
\newcommand{\B}{\mbox{\bf B}}
\newcommand{\bits}[1]{(#1)}
\newcommand{\GW}{\mathop{\mbox{\sc GW}}\nolimits}

\newcommand{\Ext}{\mathop{\mbox{\sc E}}\nolimits}
\newcommand{\Disp}{\mathop{\mbox{\sc D}}\nolimits}
\newcommand{\Con}{\mathop{\mbox{\sc C}}\nolimits}
\newcommand{\Samp}{\mathop{\mbox{\sc S}}\nolimits}

\renewcommand{\H}{{\rm H}}
\newcommand{\Hash}{{{\cal H}}}
\newcommand{\GF}{{{\rm GF}}}
\newcommand{\AP}[1]{{\rm AP}_{#1}}
\newcommand{\LD}{{\rm LD}}

  \title{Entropy waves, the zig-zag graph product,\\ and new constant-degree
  expanders}

\shorttitle{The zig-zag graph product} 

 \acknowledgements{Part of the research of the first author was performed while  
visiting the Institute for Advanced Study, Princeton, NJ.  Work of the second author was done while at MIT, supported by
an {NSF} Mathematical Sciences Postdoctoral Research Fellowship.  The third author was partially supported by NSF
grants CCR-9987007 and CCR-9987845. \hfill\break
\hglue23pt {\it Key words}\/: expander graphs, graph products, entropy.}
 \twoauthors{Omer Reingold, Salil Vadhan,}{Avi Wigderson}
 \institutions{ AT\&T Labs - Research, Florham Park, NJ\\
{\eightpoint {\it E-mail address\/}: omer@research.att.com}\\
\vglue6pt
Division of Engineering \& Applied Sciences, Harvard University, Cambridge, MA\\
{\eightpoint {\it E-mail address\/}: salil@eecs.harvard.edu}\\ \vglue6pt
Institute for Advanced Study, Princeton
and The Hebrew University, Jerusalem, Isreal\\
{\eightpoint {\it E-mail address\/}: avi@math.ias.edu}}
 
  \centerline{\bf Abstract}
\vglue6pt
The main contribution
 of this work is a new type of graph product, which we
call the {\it zig-zag product}.  
Taking a product of a large graph with a small graph, the resulting
graph inherits (roughly) its size from the large one, its degree from the
small one, and its expansion properties from both! Iteration
yields simple explicit constructions of constant-degree
expanders of arbitrary size, starting from one constant-size expander.

Crucial to our intuition (and simple analysis) of the properties of
this graph product is the view of expanders as functions which act as
``entropy wave'' propagators --- they transform probability
distributions in which entropy is concentrated in one area to
distributions where that concentration is dissipated. In these terms,
the graph product affords the constructive interference of two such
waves.

Subsequent work~\cite{ALW01}, \cite{MW01}
relates the zig-zag product of graphs to the standard
semidirect product of groups, leading to new results and
constructions on expanding Cayley graphs.

\bigbreak \centerline{\bf Contents}
\begin{small}
\sni{1} Introduction
\ssni{1.1} Expander graphs
\ssni{1.2} Overview of expander construction
\ssni{1.3} The zig-zag graph product
\ssni{1.4} Intuition
\ssni{1.5} Expanders and extractors
\ssni{1.6} Extensions to the expander construction
\ssni{1.7} Subsequent work: Connections with semidirect product in groups
\ssni{1.8} Organization of the paper

\sni{2} {Preliminaries}
\ssni{2.1} Graphs and rotations
\ssni{2.2} Eigenvalues and expansions
\ssni{2.3} Squaring and tensoring

\sni{3} The zig-zag product and the expander construction
\ssni{3.1} The zig-zag graph product
\ssni{3.2} The recursion
\sni{4} Analysis of the zig-zag product
\ssni{4.1} The basic eigenvalue bound
\ssni{4.2} Improved analysis of the eigenvalue 

\sni{5} The base graph
\ssni{5.1} The affine plane
\ssni{5.2} Low-degree polynomials

\sni{6} Variants on the zig-zag theme

\ssni{6.1} A ``derandomized" zig-zag product
\ssni{6.2} The replacement product
\end{small}

 \section{Introduction} \label{intro:sec}

 1.1. {\it Expander graphs}.
Expanders are graphs which are sparse but nevertheless highly connected.
A precise definition will be given in the next section, 
but here we informally list some
properties of such graphs (which are equivalent when
formally stated and can serve as alternate definitions)
\begin{itemize}
\item The graph satisfies ``strong'' isoperimetric inequalities.
\item Every set of vertices has ``many'' neighbors.
\item Every cut has ``many'' edges crossing it.
\item A random walk on the graph 
converges quickly to the stationary distribution.
\end{itemize}

Expander graphs have been used to address many fundamental
problems in computer science, on topics including
network
design (e.g. \cite{Pip87}, \cite{PY82}, \cite{AKS83}),
complexity theory (\cite{Val77}, \cite{Sip88}, \cite{Urq87}),
derandomization (\cite{NN93}, \cite{INW94}, \cite{IW97}),
coding theory (\cite{SS96}, \cite{Spi96}), and cryptography
(\cite{GILVZ90}). Expander graphs have also found some
applications in various areas of pure
mathematics~\cite{KR83}, \cite{Lub94}, \cite{Gro00}, \cite{LP01}.

Standard probabilistic arguments (\cite{Pin73})
show that almost every
constant-degree ($\geq 3$) graph is an expander.
However, explicit and efficient
construction of such graphs (which is required by most of the computer
science applications
above) seems to be much harder. This problem leads to
an exciting and extensive body of research,
developed mainly by mathematicians
intrigued by this computer science challenge.

Most of this work was guided by the 
algebraic characterization of expanders, developed in
\cite{Tan84}, \cite{AM85}, \cite{Alo86a}. They showed the intimate relation of
(appropriate quantitative versions of) 
all the properties above to the spectral gap in the adjacency matrix
(or, almost equivalently, the Laplacian) of the graph.
Using it, expanders can be defined as follows:
An infinite family $G_n$ of $D$-regular graphs is an {\it expander family}
if for all $n$ 
the second largest\break (in absolute value) eigenvalue of the adjacency
matrix of $G_n$ is bounded {\it uniformly} 
from above by the same $\lambda <D$. (Note that the degree $D$ is
independent of $n$; this is what we mean by ``constant degree.'')\footnote{On an intuitive level, the connection between
the spectral
  gap and the combinatorial and probabilistic properties of expanders
  listed above should not be surprising. For example, it is well known
that the standard random walk on the graph converges exponentially 
with base $\lambda /D$ to the stationary uniform
distribution. Moreover, equal partitions of the vertices of a graph,
thought of as $\pm1$-vectors, are orthogonal to the uniform distribution, and
so the bilinear form representing the number of edges in the cut can
be bounded in terms of the gap between $D$ and $\lambda$.}

This algebraic definition
naturally led researchers to consider algebraic constructions, where
this eigenvalue can be estimated. The celebrated sequence of papers
\cite{Mar73},  \cite{GG81},  \cite{AM85},  \cite{AGM87},  \cite{JM87},  
\cite{LPS88},  \cite{Mar88},\break  \cite{Mor94}  
provided such constant-degree expanders. 
All these graphs are very simple to describe: given the name of a
vertex (in binary), its neighbors can be computed in polynomial time
(or even logarithmic space). This level of explicitness is essential for many
of the applications. However, the
analysis bounding the eigenvalue
is quite sophisticated (and often based on deep mathematical
results). Thus, it is hard to intuitively understand why these graphs
are expanders.

A deviation from this path was taken in \cite{Ajt94},
where a
combinatorial construction of cubic expanders was proposed. It starts with an
arbitrary cubic $N$-vertex graph and applies a sequence of polynomially many
local operations which gradually increase the girth and turn it into
an expander. However, the resulting graphs do not have any simply described
form, and they lack the explicitness level (and hence applicability)
of the algebraic constructions mentioned above.

In this work, we give a simple, combinatorial construction of
constant-degree expander graphs.
Moreover, the analysis proving
expansion (via the second eigenvalue) is as simple and follows a
clear intuition. The construction is iterative,
and needs as a basic
building block a {\it single, almost arbitrary} expander of constant
size. The parameters required from it can be easily obtained
explicitly, but exhaustive search
is an equally good solution since it requires only constant time.
Simple operations applied to this graph generate another whose size is
increased but whose degree and expansion remain unchanged. This process
continues, yielding arbitrarily large expanders. 

The heart of the
iteration is our new ``zig-zag'' graph product.
Informally, taking a product of a large graph with 
a small graph, the resulting
graph inherits (roughly) its size from the large one, its degree from the
small one, and its expansion properties from both!  (That is, the
composed graph has good expansion properties as long as the two original
graphs have good expansion properties.)

In the next subsections we give high level descriptions of the
iterative construction, the new graph product, the intuition behind it,
various  extensions. We then mention subsequent work on the relation of
the zig-zag product in graphs to the semidirect product in groups and its
applications to expanding Cayley graphs.
 
\phantom{WEIRD}
 1.2. {\it Overview of expander construction}.
In this section, we describe a simplified, but
less efficient, version of our expander construction and omit
formal proofs.   Our full construction is described in detail
in Section~\ref{expand-main:sec}.
Throughout this section, all graphs are regular, 
undirected, and may have loops and
parallel edges. The {\it adjacency matrix} of an $N$-vertex
graph $G$ is the matrix
$M$ whose $(u,v)^{\rm th}$ entry
is the number of edges between vertices $u$ and $v$. If the graph
is $D$-regular, then  the {\it normalized adjacency matrix} 
is simply $M/D$. Note
that this stochastic matrix is the transition probability matrix of
the natural random walk on $G$, every step of which moves a ``token''
from a current vertex along a uniformly chosen edge to a  neighboring
vertex. It is easy to see that this matrix
has an eigenvalue of 1, corresponding to the constant eigenvector,
and it turns out that all other 
eigenvalues have absolute value less than
1.  Our primary interest will be
the second largest
(in absolute value) eigenvalue (which is known to govern the
convergence rate of the random walk, and as mentioned above is the
essence of expansion).

Thus, three essential parameters play a role in an expander --- size, degree
and expansion. We classify graphs accordingly.

\numbereddemo{Definition}  
An {\it $(N,D,\lambda )$-graph} is any $D$-regular graph on $N$ vertices,
whose normalized adjacency matrix has 
second largest (in absolute value) eigenvalue at most $\lambda$.
\enddemo

 {The basic operations.}
We use two operations on (the adjacency matrices of)
graphs --- the standard matrix squaring, and our new zig-zag
graph product. Here is their effect on these three parameters.

 \vglue12pt {\sc Squaring:} Let $G^2$ denote the square of $G$. Then

 \numbereddemo{Fact}  $(n,d,\la )^2 \ra (n,d^2,\la ^2)$.
 \enddemo

{\sc The zig-zag product:}  Let $G_1  \zigzag G_2$
denote the new graph product
  of $G_1$ and $G_2$. Then
 \pagebreak

\proclaim{Theorem}\label{zigzag-intro}
$(N_1,D_1, \lambda_1 ) \zigzag (D_1, D_2, \lambda_2 ) \ra
(N_1\cdot D_1,\, D_2^2,\, \lambda_1+ \lambda_2+\lambda_2^2).$
\endproclaim
 
 \vglue-12pt
(The eigenvalue bound of $\lambda_1+\lambda_2+\lambda_2^2$ is
improved somewhat in Sections~\ref{expand-main:sec} and 
\ref{improve-evalue:sec}.)

\medbreak {\it The iterations.}
Let $H$ be any $(D^4,D,1/5)$-graph, which will serve as the
building block for our construction. We define
a sequence of graphs $G_i$ as follows.
\begin{itemize}
\item[$\bullet$] $G_1=H^2$,

\item[$\bullet$]  $G_{i+1}= G_i^2 \zigzag H$.
\end{itemize}
{From} Fact~1.2 and Theorem~1.3
above, it is easy to conclude that 
this sequence is indeed an infinite family of expanders:
\proclaim{Theorem}
For every $i${\rm ,} $G_i$ is an $(N_i, D^2, 2/5)$\/{\rm -}\/graph with $N_i=D^{4i}$.
\endproclaim  

 \vglue-12pt
This construction is not as efficient as we would like
--- computing neighborhoods in $G_i$ takes time polynomial
in $N_i$ 
rather than polynomial in $\log N_i$. 
 As we show in Section~\ref{expand-main:sec},
this is easily overcome by augmenting the iterations with another standard
graph operation.

\medbreak 1.3. {\it The zig-zag graph product}.
The new product mentioned above takes a large graph and a small one, and
produces a graph that (roughly speaking) inherits
the size of the large one but the
degree of the small one. This was the key to creating arbitrarily
large graphs  with bounded degrees. Naturally, we are concerned with
maintaining the expansion properties of the two graphs. First, we
describe the product.

For simplicity, we assume
that the edges in our $D$-regular graphs
are\break $D$-colored; that is, they are
partitioned to $D$ perfect matchings.
(This assumption loses generality, and we will remove it in the formal
construction in \S \ref{expand-prelim:sec}.) For
a color $i\in[D]$ and a vertex $v$ let $v[i]$ be the neighbor of $v$
along the edge colored~$i$.
With this simple notation, we can formally define the zig-zag
product $\zigzag$ (and then explain it).

\numbereddemo{Definition}
Let $G_1$ be an $D_1$-regular graph on $[N_1]$ and $G_2$ a\break $D_2$-regular graph
on $[D_1]$. Then $G_1\zigzag G_2$ is a $D_2^2$-regular graph on
$[N_1]\times [D_1]$
defined as follows: For all $v\in [N_1], k\in [D_1], i,j\in [D_2]$,
the edge $(i,j)$ connects the vertex $(v,k)$ to the vertex
$(v[k[i]],k[i][j])$.
\enddemo

What is going on?
Note that the size of the small graph $G_2$ is the degree of the
large graph $G_1$. Thus a vertex name in $G_1\zigzag G_2$
has a first component
which is a vertex of the large graph, and a second which is viewed
both as a vertex of the small graph {\it and} an edge color of the large
one. The edge label in $G_1\zigzag G_2$ is just a pair of edge labels in
the small graph.
One step in the new product graph from a vertex $(v,k)$
along the edge $(i,j)$ can be broken into three substeps.
\begin{itemize}
\item[1.] $(v, k) \ra (v,k[i])$  --- A step (``zig'') in the small graph moving
$k$
  to $k[i]$. This affects only the second component, according to the
  first edge label.
\item[2.] $(v,k[i]) \ra (v[k[i]],k[i])$ --- A step in the large graph,
  changing the first component according to the second, viewed as an
  edge color.
\item[3.] $(v[k[i]],k[i]) \ra (v[k[i]],k[i][j])$ -- A step (``zag'') in the small
graph moving $k[i]$
  to $k[i][j]$. This affects only the second component, according to the
  second edge label.
\end{itemize}

\demo{{\rm 1.4.} Intuition}
Why does it work? More precisely, why does Theorem
1.3  hold?  What this theorem says intuitively, is
that $G_1\zigzag G_2$ is a good expander as long as both $G_1$ and $G_2$ are
good expanders.
Consider the above three steps as a random walk on $G_1\zigzag G_2$.
Then Steps 1 and 3 are independent random  steps on the small graph.
If at least one of them ``works'' as well as it does in the small
graph, then this would guarantee that the new graph is as good an expander as
the small one. So let us argue (very
intuitively) that indeed one of them ``works.''

A random step in an expander increases the ($\H_2$-) entropy of a
distribution on the vertices, {\it provided that it is not already too close
to uniform.}
Let us consider a distribution on the vertices of the new
graph $(v,k)$. Roughly speaking, there are two cases.
\begin{itemize}
\item If the distribution of the second component $k$ 
(conditioned on $v$) is not too uniform,
then Step~1 ``works.''
Since Step~2 is just a permutation and Step~3 is a random step
on a regular graph, these steps cannot make the distribution
less uniform and undo the progress made in Step~1.
\item
If $k$ (conditioned on $v$) is very close to uniform, then
Step~1 is
a ``waste.'' However, Step~2 is then like a real random step in the
large expander $G_1$! This means that the entropy of the first component $v$
increases.
Note that Step~2 is a permutation on
the vertices of $G_1\zigzag G_2$, so if entropy increases in the first
component, it decreases in the second.
That means that in Step~3 we are in the good case (the conditional
distribution on the second component is far from uniform), and the
entropy of the second component will increase by the expansion of the
small graph.
\end{itemize}

\vglue6pt 
The key to this product is that Step~2 is simultaneously
a permutation (so that any progress made in Step~1 is preserved)
and an operation whose ``projection'' to the first component
is simply a random step on the large graph
 (when the second component
is random).
All previous discussions of expanders focused on the
increase of entropy to the vertex distribution by a step along a
random edge. We insist on keeping track of that edge name, and
consider the joint distribution! In a good expander, if the edge is
indeed random, the entropy propagates from it to the vertex. This
reduces the (conditional) entropy in the edge. Thus the ``entropy wave''
in Step~2, in which
no fresh randomness enters the distribution on vertices
of $G_1\zigzag G_2$, is what
facilitates  entropy increase in Steps 1 or 3.
Either the ``zig'' step does it, if there is room for more entropy in
$k$, or if not (which may be viewed as destructive interference
of the large and small waves in
Step~1), Step~2 guarantees constructive interference in Step~3. 
Moreover, Step~1 is not redundant as, if there is no or little initial entropy in
$k$, the wave of Step~2 (being a permutation) may flood $k$ with
entropy, destroying the effect of Step~3.

The formal proof of Theorem 1.3 follows this
intuition quite closely, and separately analyzes these two extreme cases.
Indeed, since it becomes linear algebra,  these two cases are very
natural to define,
and the only ones to worry about ---
all intermediate cases follow by linearity! 
Moreover, the variational definition of the second eigenvalue better
captures the symmetry of the zig and zag steps (and gives a better
bound than what can be obtained from this asymmetric intuition).

\vglue6pt

\demo{{\rm 1.5.} Expanders and extractors}
Here we attempt an intuitive explanation of how we stumbled on the definition
of the zig-zag product, and the intuition that it does what it should.
While this subsection may not be self-contained, it will at least lead
the interested reader to discover more of the fascinating world of
extractors.

The current paper is part of research described in our conference paper
\cite{RVW00}
which deals with constructions of both expanders and extractors.
Extractors are combinatorial objects, defined by \cite{NZ96},
which, roughly speaking, ``purify''
arbitrary nonuniform probability distributions into uniform ones.
These objects  are as fascinating and as applicable as expanders 
(see, e.g., the  survey papers \cite{Nis96}, \cite{NT99}).
Like expanders, their applications demand explicit construction.
Like with expanders, the quest for such constructions has been
extremely fruitful and illuminating for complexity theory. Unlike
expanders, the construction of optimal extractors is still a
challenge, although the best existing ones are quite close to optimal
(see the current state of the art, as well as a survey of previous
constructions, in \cite{RSW00}, \cite{TUZ01}).

Expander graphs were ingredients in some previous
extractor constructions (as
extractors may be viewed as graphs as well). Here the situation is reversed.
The expander construction of this paper 
{\it followed} our discovery of nearly optimal {\it high
  min-entropy} extractors, which handle the ``purification'' of
  distributions which are already not too far from being  uniform.
A key idea in approaching optimality (following \cite{RR99}) was
preserving the unused entropy in a random step on an extractor. 
This lead to a (more complex) type of zig-zag product, and from
it, iterative constructions of such extractors.
Translating this idea to the expander world turned out to be cleaner
and more natural than in the extractor world. It led to our
understanding of the role of the edge-name as a keeper of the unused
entropy in a step of a 
standard random walk, and to the zig-zag product defined above. \enddemo

  1.6. {\it Extensions to the expander construction}.
The list below details the extensions and refinements we
obtain to the basic expander construction outlined above. All these
will be part of the formal sections which follow.

\vglue6pt {\it More explicit graphs.}
As mentioned above, this construction is not as efficient as we would like
--- computing neighborhoods in $G_i$ takes time 
polynomial in $N_i$ rather than in $\log N_i$.   As we show in Section~\ref{expand-main:sec},
this is easily overcome by augmenting the iterations with another standard
graph operation, namely taking tensor powers of the adjacency matrix.

\vglue6pt {\it Describing graphs by {\rm ``}\/rotation maps}.''
Another explicitness problem in the simple construction above is the
assumption that the our $D$-regular graphs are given together with a
proper $D$-coloring of the edges. This property is not preserved by
the zig-zag product. To avoid it, we describe graphs more generally by
their ``rotation maps,'' and show how this description is explicitly
preserved by all graph operations in our construction.

\vglue6pt {\it Smaller degree.}
  A naive and direct implementation of our graph
  product yields expanders whose degree is reasonable, but not that
  small (something under 1000). 
  In Section~\ref{recursion:sec}, we show how to combine this
  construction, together with {\it one, constant-size} 
        cycle,   
  to
  obtain an infinite family of explicit degree 4 expanders. Again,
  this combination uses the zig-zag product.  In fact, using
  the replacement product described below, we obtain explicit
  degree 3 expanders (which is the smallest possible).

\vglue6pt {\it Choice of the base graph.}
Our expander construction requires an initial ``constant size''
base graph $H$ as a building block.  While exhaustive search can
be used to find such an $H$ (since it is constant size), for completeness we
include two elementary explicit
constructions (from \cite{Alo86b,AR94})
which can be used instead.

\vglue6pt {\it Better degree vs.\ eigenvalue relation.}
The best relationship between degree and second largest eigenvalue is
obtained by {\it Ramanujan} graphs, in which the second eigenvalue is
$2\sqrt{D-1}/D$.  This equals the first eigenvalue of the\break
$D$-regular infinite tree, and it 
is known
that no  finite $D$-regular graph
can have a smaller second largest eigenvalue
(cf., \cite{Alo86a}, \cite{LPS88}, \cite{Nil91}).
Remarkable graphs achieving this optimal bound
were first constructed independently 
by \cite{LPS88} (who coined the term
Ramanujan graphs) and by \cite{Mar88}.
 \pagebreak

Our constructions do not achieve this tight relationship.
The  zig-zag product, applied recursively to one fixed Ramanujan graph,
  will yield $D$-regular expanders of second largest 
  eigenvalue $O(1/D^{1/4})$. A ``partially derandomized'' variant
  of our zig-zag product, given in 
        Section~\ref{zigzagmod:sec}, improves this relation and
  achieves second eigenvalue $O(1/D^{1/3})$.

\vglue2pt {\it A simpler product.}
Perhaps the most natural way to combine $G_1$ with $G_2$ when the size
of $G_2$ is the degree of $G_1$ is simply replace every vertex of
$G_1$ with a copy of $G_2$ in the natural way, keeping the edges of
both graphs. This {\it replacement product}, which was often used for
degree-reduction purposes (e.g., when $G_2$ is a cycle the resulting
graph has degree 3) turns out to enjoy similar properties of the
zig-zag product: if both $G_1$ and $G_2$ are expanders, then so is their
replacement product. Moreover, the proof is by a reduction --- the
zig-zag product is a subgraph of the cube ($3^{\rm rd}$ power) of the
replacement product, immediately giving an eigenvalue bound.

\demo{{\rm 1.7.} Subsequent work\/{\rm :} Connections with semidirect product in groups}
Subsequent to this work, it was shown in \cite{ALW01} that
the zig-zag (and replacement) products can be viewed as a
generalization of the standard semidirect product of groups. This was
used in \cite{ALW01} to construct a family of groups which is
expanding with one (constant size) set of generators, but is not
expanding with another such set. The connection was further developed
in \cite{MW01} 
to produce new families of expanding Cayley
graphs, via bounds on the number of irreducible representations of
different dimensions in terms of the expansion.
\enddemo

1.8. {\it Organization of the paper}.
In Section~\ref{expand-prelim:sec}, we give preliminary definitions and
basic facts. In Section~\ref{expand-main:sec}, we define the zig-zag
graph product, describe the construction of expanders, and state their
properties. In particular, it deals with the first four
``extensions'' listed in the previous subsection. In
Section~\ref{analysis:sec}, we analyze the expansion of the zig-zag
product. 
In Section~\ref{basegraph:sec}, we discuss some ways to obtain
the base graph used in our expander construction.
In Section~\ref{zigzagmod:sec}, we give two extensions to the
basic zig-zag product. The first
is a ``derandomized''  variant of our
basic zig-zag product, which enjoys a better relationship between the
degree and the expansion. The second is the simple, natural
{\it replacement} product. 

\vglue-12pt
\section{Preliminaries} \label{expand-prelim:sec}
\vglue-6pt

2.1. {\it  Graphs and rotations}.
All graphs we discuss may have self-loops and parallel edges. They are
best described by their (nonnegative, integral) 
adjacency matrix. Such a graph is
{\it undirected} if and only if the adjacency matrix is symmetric. It is 
{\it $D$-regular} if
the sum of entries in each row (and column) is $D$ (so exactly $D$
edges are incident to every vertex).

Let $G$ be a $D$-regular undirected graph on $N$ vertices.  
Suppose that
the edges leaving each vertex of $G$ are labeled from $1$ to $D$ in
some arbitrary, but fixed, way.  Then for $v,w\in [N]$ and
$i\in [D]$, it
makes sense (and is standard) to say
``the $i^{\rm th}$ neighbor of vertex $v$ is $w$.''
In this work, we make a point to always keep track of the
edge traversed to get from $v$ to $w$.  This is formalized as follows:
\numbereddemo{Definition}\label{rotation:def}
For a $D$-regular undirected graph $G$, the
{\it rotation map} $\Rot{G} : [N]\times [D] \to [N]\times
[D]$
is defined as follows:
$\Rot{G}(v,i)=(w,j)$ if
the $i^{\rm th}$ edge incident to $v$ leads to $w$,
and this edge is the $j^{\rm th}$ edge incident to $w$.
\enddemo

This definition enables us to remove the simplifying assumption made
in the introduction, which was 
that the label of an edge is the same from the perspective
of both endpoints, i.e. $\Rot{G}(v,i)=(w,j) \Rightarrow i=j$.
{From} Definition~\ref{rotation:def}, 
it is clear that $\Rot{G}$ is a permutation, and
moreover $\Rot{G}\circ \Rot{G}$ is the identity map.  

We will always
view graphs as being specified by their rotation maps.  Hence
we call a family
$\GraphFam$ of graphs {\it explicit}
if for every $G\in \GraphFam$,
$\Rot{G}$ is computable in time $\poly(\log N)$, where $N$ is
the number of vertices of $G$.
That is, graphs in $\GraphFam$ are indexed by some parameters
(such as the number of vertices and the degree, which may
be required to satisfy some additional relations) and there
should be a single algorithm which efficiently computes
$\Rot{G}$ for any $G\in\GraphFam$
when given these parameters as an additional input. The notation
$\poly()$ stands for a fixed (but unspecified) 
polynomial function in the given variables.
We will often informally refer to an individual graph as explicit,
as shorthand for saying that the graph comes from an explicit family.

Our constructions will be iterative (or recursive), and will be based
on a sequence of composition operations,
constructing new graphs from given ones. 
The definition of these compositions (or products) will show
how the rotation map
of the new graph can be computed using ``oracle access'' to the rotation
maps of the given graphs. (By giving an algorithm ``oracle access'' to
a function~$f$, we mean that the algorithm is given power
to evaluate $f$ on inputs of its choice at the cost
of one time step per evaluation.)
Given the time complexity of such a computation
{\it and} the number of oracle calls made, it will be easy to compute
the total time required by a recursive construction. 

\demo{{\rm 2.2.} Eigenvalues and expansions}
The {\it normalized adjacency matrix} $M$ of $G$
is the adjacency matrix of $G$ divided by $D$.
In terms of the rotation map, we have:
$$M_{u,v}=\frac{1}{D}\cdot
\left|\{ (i,j)\in [D]^2 : \Rot{G}(u,i)=(v,j)\}\right|.$$
$M$ is simply the transition matrix of a random walk on $G$.
By the $D$-regularity of $G$, the all-1's vector
$\NOnes=(1,1,\ldots,1)\in \R^N$ is an eigenvector of $M$ of
eigenvalue 1.  It is turns out that all the other eigenvalues
of $M$ have absolute value at most 1, and it is
well known that the second largest eigenvalue of $G$
is a good measure of $G$'s expansion
properties~\cite{Tan84}, \cite{AM85}, \cite{Alo86a}.
We
will use the following variational characterization of the second
largest eigenvalue.
\enddemo

\numbereddemo{Definition}\label{evalue:def}
$\lambda(G)$ denotes the {\it second largest eigenvalue}
(in absolute value) of $G$'s
normalized adjacency matrix.  Equivalently,
$$\lambda(G) =
\max_{\ALPHA \perp \NOnes}
\frac{|\langle \ALPHA, M\ALPHA\rangle|}{\langle \ALPHA, \ALPHA\rangle}
= \max_{\ALPHA \perp \NOnes}
\frac{\|M\ALPHA\|}{\|\ALPHA\|}.$$
\enddemo

Above, $\langle\cdot,\cdot\rangle$ refers to the standard inner
product in $\R^N$ and $\|\ALPHA\|=\sqrt{\langle \ALPHA,\ALPHA\rangle}$.

The meaning of $\lambda(G)$ can be understood as follows:
Suppose $\PI\in[0,1]^N$ is a probability distribution on the vertices
of $G$.  By linear algebra, 
$\PI$ can be decomposed as $\PI=\NUni+\PIperp$, where
$\NUni=\NOnes/N$ is the uniform distribution and $\PIperp\perp\NUni$.
Then $M\PI=\NUni+M\PIperp$ is the probability distribution on vertices
obtained
by selecting a vertex $v$ according to $\PI$ and then moving to a
uniformly
selected neighbor of $v$.  By Definition~\ref{evalue:def},
$\|M\PIperp\|\leq \lambda(G)\cdot \|\PIperp\|$.  Thus $\lambda(G)$
is a measure of how quickly the random walk on $G$ converges to the
uniform distribution.
Intuitively, the smaller $\lambda(G)$ is, the better the expansion
properties of
$G$.
Accordingly, an (infinite) family $\GraphFam$ of graphs is called
a family of {\it expanders}
if these eigenvalues are bounded away from 1, i.e.
there is a constant $\lambda<1$ such that $\lambda(G)\leq \lambda$
for all $G\in \GraphFam$.  It was shown by
Tanner~\cite{Tan84} and Alon and Milman~\cite{AM85}
that this implies (and is in fact equivalent to~\cite{Alo86a})
the standard notion
of {\it vertex expansion}: there is a constant
$\varepsilon>0$ such that
for every $G\in \GraphFam$ and for any set $S$ of at most
half the vertices
in $G$, at least $(1+\varepsilon)\cdot |S|$
vertices of $G$ are connected to some
vertex in $S$.

As mentioned in the introduction,
we refer
to a $D$-regular undirected graph $G$
on $N$ vertices such that $\lambda(G)\leq \lambda$ as an
{\it $(N,D,\lambda)$-graph}.
Clearly, achieving expansion is easier as the degree gets larger.  The
main
goal in constructing expanders is to minimize the degree, and, more
generally, obtain the best degree-expansion tradeoff.  Using
the probabilistic method, Pinsker~\cite{Pin73} showed that
most 3-regular graphs are expanders (in the sense of vertex expansion),
and this result was extended to eigenvalue bounds
in \cite{Alo86a}, \cite{BS87}, \cite{FKS89}, \cite{Fri91}.
The best-known bound on the eigenvalues of random graphs
is due to Friedman~\cite{Fri91}, who showed that
most $D$-regular graphs have second largest eigenvalue
at most $2/\sqrt{D}+O((\log D)/D)$ (for even $D$).
In fact, the bound of $2\sqrt{D-1}/D$ is the best
possible
for an infinite family of graphs,
as shown by Alon and Boppana (cf.\ \cite{Alo86a}, \cite{LPS88}, \cite{Nil91}).
Graphs whose second largest eigenvalue meets this optimal bound are
called {\it Ramanujan graphs}. It is easy to verify that this
value is the {\it largest} eigenvalue of the random walk on the {\it infinite} $D$-regular tree.

While these probabilistic arguments provide strong existential results,
applications of expanders in computer science often require
{\it explicit} families of constant-degree expanders.
The first
such construction was given by\break Margulis~\cite{Mar73}, with
improvements and simplifications by Gabber and\break
Galil~\cite{GG81}, Jimbo and Maruoka~\cite{JM87},
Alon and Milman~\cite{AM85},
and Alon, Galil, and Milman~\cite{AGM87}.  Explicit 
families of Ramanujan graphs
were first constructed by
Lubotzky, Phillips, and Sarnak~\cite{LPS88} and
Margulis~\cite{Mar88}, with more recent constructions
given by Morgenstern~\cite{Mor94}. 
The best eigenvalues
we know how to achieve using our approach 
are $O(1/D^{1/3})$.

\demo{{\rm 2.3.} Squaring and tensoring}
In addition to the new zig-zag product, our 
expander construction makes use of two standard operations on
graphs --- squaring and tensoring.  Here we describe these
operations in terms of rotation maps and state their effects on
the eigenvalues.

Let $G$ be a $D$-regular multigraph on $[N]$ given by
rotation map $\Rot{G}$.  The {\it $t^{\rm th}$
power} of $G$ is the $D^t$-regular graph $G^t$ whose rotation
map is given by $\Rot{G^t}(v_0,(k_1,k_2,\ldots,k_t))=
(v_t,(\ell_t,\ell_{t-1},\ldots,\ell_1))$, where
these values are computed via the rule
$(v_i,\ell_i) = \Rot{G}(v_{i-1},k_i)$.\enddemo

\proclaim{Proposition}\label{power}
If $G$ is an $(N,D,\lambda)${\rm -}graph{\rm ,} then $G^t$ is an
$(N,D^t,\lambda^t)$\/{\rm -}\/graph{\rm .}  Moreover{\rm ,} $\Rot{G^t}$ is computable
in time $\poly(\log N,\log D,t)$ with $t$ oracle queries to $\Rot{G}${\rm .}
\endproclaim

\vglue-4pt
{\it Proof}.
The normalized adjacency matrix of $G^t$ is the
$t^{\rm th}$ power of the normalized adjacency matrix of $G$,
so all the eigenvalues also get raised to the $t^{\rm th}$ power.
\hfill\qed\vglue7pt

Let $G_1$ be a $D_1$-regular multigraph on $[N_1]$ and
let $G_2$ be a $D_2$-regular multigraph on $[N_2]$.
Define the {\it tensor product} $G_1\otimes G_2$ to
be the $D_1\cdot D_2$-regular multigraph on $[N_1]\times [N_2]$
given by $\Rot{G_1\otimes G_2}((v,w),(i,j))=((v',w'),(i',j'))$,
where $(v',i')=\Rot{G_1}(v,i)$ and $(w',j')=\Rot{G_2}(w,j)$.
In order to analyze this construction (and our new graph
product), we need some concepts from linear algebra.  For
vectors $\ALPHA\in\R^{N_1}$ and $\BETA\in\R^{N_2}$, their
{\it tensor product} is the vector
$\ALPHA\otimes\BETA\in\R^{N_1\cdot N_2}$ whose $(i,j)^{\rm th}$ entry
is $\ALPHA_i\cdot\BETA_j$.  If $A$ is an $N_1\times N_1$ matrix
and $B$ is an $N_2\times N_2$ matrix, then there
is a unique $N_1N_2\times N_1N_2$ 
matrix $A\otimes B$ (again called the {\it tensor
product})
such that $(A\otimes B)(\ALPHA\otimes \BETA)=(A\ALPHA)\otimes(B\BETA)$
for
all $\ALPHA,\BETA$.

\proclaim{Proposition} \label{tensor:prop} \hskip-8pt
If $G_1$ is an $(N_1,D_1,\lambda_1)$\/{\rm -}\/graph
and $G_2$ is an $(N_2,D_2,\lambda_2)$\/{\rm -}\/graph{\rm ,}
then $G_1\otimes G_2$ is an
$(N_1\cdot N_2,D_1\cdot D_2,\max(\lambda_1,\lambda_2))$\/{\rm }\/-graph{\rm .}
Moreover{\rm ,} $\Rot{G_1\otimes G_2}$ is computable
in time $\poly(\log N_1N_2,\log D_1D_2)$ with
one oracle query to $\Rot{G_1}$ and one oracle query to $\Rot{G_2}${\rm .}
\endproclaim

\demo{Proof}
The normalized adjacency matrix of $G_1\otimes G_2$ is
the tensor product of the normalized adjacency matrices
of $G_1$ and $G_2$.  Hence its eigenvalues are the pairwise
products of eigenvalues of $G_1$ and $G_2$.  The
largest eigenvalue is $1\cdot 1$, and the second
largest eigenvalue is either $1\cdot \lambda_2$ or $\lambda_1\cdot 1$.
\enddemo

\section{The zig-zag product and the expander construction} 
\label{expand-main:sec}

In the introduction, we described how to obtain a
family of expanders by iterating 
two operations on graphs --- squaring and
the new ``zig-zag'' product.
That description used a simplifying
assumption about the edge labeling.  In terms of rotation maps, the
assumption was that $\Rot{}(v,i)=(w,j)\Rightarrow i=j$.
In this section, we describe the construction in terms of arbitrary
rotation maps
and prove its properties.
The expander construction given here
will also use tensoring to
improve the efficiency to polylogarithmic
in the number of vertices. This deals with the first
two items in the ``extensions'' subsection of the introduction, which
are summarized in Theorem~\ref{zigzag:thm}. The third item,
obtaining
expanders of degree $4$, will follow in Corollary~\ref{deg4:cor}.
The analysis of the zig-zag product is deferred
to the following section.

\demo{{\rm 3.1.} The zig-zag graph product}
We begin by describing the new graph product in terms
of rotation maps.
Let $G_1$ be a $D_1$-regular multigraph on $[N_1]$ and $G_2$ a
$D_2$-regular multigraph on $[D_1]$.
Their {\it zig-zag product} is a $D_2^2$-regular
multigraph $G_1\zigzag\, G_2$ on $[N_1]\times [D_1]$.  We view every
vertex $v$
of $G_1$ as being blown up to a ``cloud'' of $D_1$ vertices
$(v,1),\ldots,(v,{D_1})$, one for each edge of $G_1$ leaving $v$.
Thus for every edge $e=(v,w)$ of $G_1$, there are two
associated vertices of $G_1\zigzag G_2$ --- $(v,k)$ and $(w,\ell)$, where
$e$ is the $k^{\rm th}$ edge leaving $v$ and the $\ell^{\rm th}$ edge leaving
$w$.  Note that these pairs satisfy the relation
$(w,\ell)=\Rot{G_1}(v,k)$.
Since $G_2$ is a graph on $[D_1]$, we can also imagine connecting
the vertices of each such cloud using the edges of $G_2$.
Now, the edges of $G_1\zigzag G_2$ are defined (informally)
as follows: we
connect two vertices $(v,k)$ and $(w,\ell)$ if it is possible to get
from $(v,k)$ to $(w,\ell)$
by a sequence of moves of the following form:
\begin{itemize}
\item[1.] Move to a neighboring vertex $(v,k')$ within
the initial cloud (using an edge of $G_2$),
\item[2.] Jump across clouds (using edge $k'$ of $G_1$) to get to
$(w,\ell')$,
\item[3.] Move to a neighboring vertex $(w,\ell)$ within
the new cloud (using an edge of $G_2$).
\end{itemize}
To make this precise, we describe how to compute the
$\Rot{G_1\zigzagsub G_2}$ given
$\Rot{G_1}$ and $\Rot{G_2}$.
\enddemo

\numbereddemo{Definition}If $G_1$ is a $D_1$-regular graph on $[N_1]$ with rotation
map $\Rot{G_1}$ and $G_2$ is a $D_2$-regular graph on $[D_1]$
with rotation map $\Rot{G_2}$, then their {\it zig-zag product}
$G_1\zigzag\, G_2$ 
is defined to be the $D_2^2$-regular graph on $[N_1]\times [D_1]$
whose rotation map $\Rot{G_1\zigzagsub G_2}$ is as follows:
\medbreak {\it $\Rot{G_1\protect\zigzagsub G_2}((v,k),(i,j))$}:
\begin{itemize}
\item[1.] Let $(k',i')=\Rot{G_2}(k,i)$,
\item[2.] Let $(w,\ell')=\Rot{G_1}(v,k')$,
\item[3.] Let $(\ell,j')=\Rot{G_2}(\ell',j)$,
\item[4.] Output $((w,\ell),(j',i'))$.
\end{itemize}
\enddemo

The important feature of this graph product is
that $G_1\zigzag G_2$ is a good expander if both $G_1$ and
$G_2$ are, as shown by the following theorem.

\proclaim{Theorem} \label{zigzag:thm}
If $G_1$ is an $(N_1,D_1,\lambda_1)$\/{\rm -}\/graph and $G_2$ is
a $(D_1,D_2,\lambda_2)$\/{\rm -}\/graph{\rm ,} then
$G_1\zigzag G_2$ is a
$(N_1\cdot D_1,D_2^2,f(\lambda_1,\lambda_2))$\/{\rm -}\/graph{\rm ,}
where
$f(\lambda_1,\lambda_2)\leq \lambda_1+\lambda_2+\lambda_2^2$ and
$f(\lambda_1,\lambda_2)<1$ when $\lambda_1,\lambda_2<1${\rm .}
Moreover{\rm ,} $\Rot{G_1\zigzagsub G_2}$ can be computed in time
$\poly(\log N,\log D_1,\log D_2)$ with one oracle query to
$\Rot{G_1}$ and two oracle queries to $\Rot{G_2}${\rm .}
\endproclaim

Stronger bounds on the function $f(\lambda_1,\lambda_2)$ are
given in Section~\ref{improve-evalue:sec}.
Before proving Theorem~\ref{zigzag:thm}, we show how it can be used to
construct an infinite family of constant-degree expanders starting from
a constant-size expander.  

\demo{{\rm 3.2.} The recursion} \label{recursion:sec}
The construction is like the construction
in the introduction, except that we
use tensoring to reduce the depth of the recursion and thereby
make the construction run in polylogarithmic time (in the size of
the graph).
 
Let $H$ be a
$(D^8,D,\lambda)$-graph for some $D$ and $\lambda$. (Various methods for
obtaining such an $H$ are described in \S \ref{basegraph:sec}.)
 For every
$t\geq 1$, we will define a $(D^{8t},D^2,\lambda_t)$-graph
$G_t$.
$G_1$ is $H^2$ and $G_2$ is $H\otimes H$.
For $t>2$, $G_t$ is recursively
defined by
$$G_t=
\left(G_{\lceil \frac{t-1}{2}\rceil} \otimes
G_{\lfloor \frac{t-1}{2}\rfloor}\right)^2 \zigzag\, H.$$

\proclaim{Theorem} \label{expfam:thm}
For every $t\geq 0${\rm ,} $G_t$ is an $(D^{8t},D^2,\lambda_t)$\/{\rm -}\/graph
with $\lambda_t=\lambda+O(\lambda^2)${\rm .}  Moreover{\rm ,}
$\Rot{G_t}$ can be computed in time $\poly(t,\log D)$
with $\poly(t)$ oracle queries to $\Rot{H}${\rm .}
\endproclaim

\demo{Proof}
A straightforward induction establishes that
the number of vertices in $G_t$ is $D^{8t}$ and that its
degree is $D^2$.
To analyze the eigenvalues, define
$\mu_t = \max\{\lambda_1,\ldots,\lambda_t\}$.
Then we have $\mu_t \leq
\max\{\mu_{t-1},\mu_{t-1}^2+\lambda+\lambda^2\}$
for all $t\geq 2$.   Solving this recurrence gives
$\mu_t\leq \lambda+O(\lambda^2)$ for all $t$.
For the efficiency, note that the depth of the recursion
is at most $\log_2 t$ and evaluating the rotation maps for
$G_t$ requires 4 evaluations of rotation maps for smaller graphs,
so the total number of recursive calls is at most $4^{\log_2 t}=t^2$.
\enddemo

In order for Theorem~\ref{expfam:thm} to guarantee
that graphs $\{G_t\}$ are expanders, 
the second largest eigenvalue $\lambda$ of the
building block $H$ 
must be sufficiently small
(say, $\lambda \leq 1/5$).  This forces
the degree of $H$ and hence the degree of the expander family
to be rather large, though still constant.
However, by zig-zagging the family $\{G_t\}$ with a cycle, we
can obtain a family of degree 4 expanders.  
More generally, we can use this method
convert any family of odd-degree expanders into a family of degree
4 expanders:

\proclaim{{C}orollary}\label{deg4:cor}
For every $\lambda<1$ and every odd $D${\rm ,} there exists a $\lambda'<1$ such that
if $G$ is an $(N,D,\lambda)$\/{\rm -}\/graph and $C$ is the cycle on $D$ vertices{\rm ,}
then $G\zigzag C$ is a $(ND,4,\lambda')$\/{\rm -}\/graph{\rm .} 
\endproclaim 

\demo{Proof}
As with any connected and nonbipartite graph, $\lambda(C)$ is strictly
less than 1 for an odd cycle $C$ 
(though $\lambda(C)\to 1$ as $D\to \infty$). Thus, the corollary
follows from Theorem~\ref{zigzag:thm}.
\enddemo

\section{Analysis of the zig-zag product}\label{analysis:sec}

This section has two subsections. In the first, we give the
basic (suboptimal)  bound of Theorem~\ref{zigzag:thm}. This bound
uses only the intuitive ideas of the introduction, and  suffices for
the construction of the previous section.
In the next, we state and prove a tighter eigenvalue bound. It uses
extra information about the  zig-zag product (which is less
intuitive). It also gives more information about the worst interplay
between the two extreme cases studied in the basic analysis, and may
hopefully shed a bit of light on the structure of the eigenvectors of
the zig-zag product.

\demo{{\rm 4.1.} The basic eigenvalue bound} \label{evalue:sec}
Now we prove Theorem~\ref{zigzag:thm}.
Recall   the intuition behind the zig-zag product.
We aim to show that for any (nonuniform)
initial probability distribution $\pi$
on the vertices of
$G_1\zigzag\, G_2$, taking a random step on $G_1\zigzag G_2$
results in a distribution that is more uniform.
We argued this intuitively in the introduction, 
by considering two extreme cases, based on the
conditional distributions induced by $\pi$ on the $N_1$ ``clouds''
of $D_1$ vertices each: one in which these conditional distributions
are far from uniform,
and the second in which they are uniform.
The actual linear algebra proof below
will restrict itself to these two cases
by decomposing any other vector into a linear combination of the two. 
Also, the argument in the introduction was not symmetric in the first and
second steps on the small graph.  Using the variational definition of
the second largest eigenvalue, we get a cleaner analysis than
by following that intuition directly.

Let $M$ be the
normalized adjacency matrix of $G_1\zigzag\, G_2$.
According to Definition~\ref{evalue:def}, we must show that,
for every vector $\ALPHA\in \R^{N_1\cdot D_1}$
such that $\ALPHA\perp \Ones{N_1D_1}$, 
$|\langle M\ALPHA,\ALPHA\rangle|$ is smaller than
$\langle \ALPHA, \ALPHA \rangle$ by a 
factor 
$f(\lambda_1,\lambda_2)$.
For intuition, 
$\ALPHA$ should be thought of as the nonuniform component
of
the probability distribution $\pi$
referred to above,  i.e.
$\pi=\Uni{N_1D_1}+\ALPHA$, where $\Uni{N_1D_1}=\Ones{N_1D_1}/N_1D_1$ is
the uniform distribution on $[N_1D_1]$.  Thus, we are showing that
$\pi$ becomes more uniform after a random step on $G_1\zigzag G_2$.

For every $v\in [N_1]$, define
$\ALPHA_v\in\R^{D_1}$ by $(\ALPHA_v)_k=\ALPHA_{vk}$.
Also define a (linear) map $\Clouds : \R^{N_1\cdot D_1} \to
\R^{N_1}$
by $(\Clouds\ALPHA)_v=\sum_{k=1}^{D_1} \ALPHA_{vk}$.
Thus,
for a probability distribution $\pi$ on the vertices
of $G_1\zigzag G_2$, $\pi_v$ is a multiple of the
conditional distribution on ``cloud $v$'' and
$\Clouds\pi$ gives the marginal distribution
on set of clouds.
By definition,
$\ALPHA=\sum_v \Unit{v}\otimes \ALPHA_v$, where
$\Unit{v}$ denotes the $v^{\rm th}$ standard basis vector in $\R^{N_1}$.
By basic linear algebra, every $\ALPHA_v$ can be
decomposed (uniquely) into $\ALPHA_v=\ALPHA_v^\parallel+\ALPHA_v^\perp$
where
$\ALPHA_v^\parallel$ is parallel to $\Ones{D_1}$ (i.e., all of its
entries
are the same) and
$\ALPHA_v^\perp$ is orthogonal to $\Ones{D_1}$ (i.e., the
sum of its entries are 0).
Thus, we obtain a decomposition of $\ALPHA$:
\begin{eqnarray*}
\ALPHA &=&
\sum_v \Unit{v}\otimes \ALPHA_v\\
&=& \sum_v \Unit{v}\otimes \ALPHA_v^\parallel
+ \sum_v \Unit{v}\otimes \ALPHA_v^\perp\\
&\eqdef& \ALPHApar+\ALPHAperp.
\end{eqnarray*}

This decomposition corresponds to the two cases in our intuition:
$\ALPHApar$
corresponds to a probability distribution
on the vertices of $G_1\zigzag G_2$ such that
the conditional distributions on the clouds
are all uniform.
$\ALPHAperp$
corresponds to a distribution
such that
the conditional distributions on the clouds
are all far from uniform.
Another way of matching
$\ALPHApar$ with the intuition is
to note that $\ALPHApar=\Clouds\ALPHA\otimes \Ones{D_1}/D_1$.
Since $\ALPHA$ and  $\ALPHAperp$ are both orthogonal to
$\Ones{N_1D_1}$,
so is $\ALPHApar$ and hence also $\Clouds\ALPHA$ is orthogonal
to $\Ones{N_1}$.

To analyze how $M$ acts on these two vectors,
we relate $M$ to the normalized adjacency matrices
of $G_1$ and $G_2$, which we denote by $A$ and $B$, respectively.
First, we decompose $M$
into the product of three matrices, corresponding to the
three steps in the definition of $G_1\zigzag G_2$'s edges.
Let $\liftB$ be the (normalized) adjacency matrix
of the graph on $[N_1]\times [D_1]$
where we connect the vertices within each cloud
according to the edges of $G_2$.
$\liftB$ is related to $B$ by the relation
$\liftB=I_{N_1}\otimes B$, where
$I_{N_1}$
is the
$N_1\times N_1$ identity matrix.
Let $\liftA$ be the permutation matrix
corresponding to
$\Rot{G_1}$.
The relationship between $\liftA$ and $A$ is somewhat
subtle, so we postpone describing it until later.
By the definition of $G_1\zigzag G_2$,
we have $M=\liftB\liftA\liftB$.
Note that both $\liftB$ and $\liftA$ are symmetric matrices, due to the
undirectedness of $G_1$ and $G_2$.

Recall that we
want to bound
$|\langle M\ALPHA,\ALPHA\rangle|/\langle \ALPHA,\ALPHA\rangle$.
By the symmetry of $\liftB$,  
\begin{equation}  \label{symmetrize:eqn}
\langle M\ALPHA,\ALPHA\rangle = 
\langle \liftB\liftA\liftB \ALPHA,\ALPHA\rangle
= \langle \liftA\liftB \ALPHA,\liftB\ALPHA\rangle.
\end{equation}
Now note that $\liftB\ALPHApar=\ALPHApar$, because
$\ALPHApar=\Clouds\ALPHA\otimes \Ones{D_1}/D_1$,
$\liftB=I_{N_1}\otimes B$, and
$B\Ones{D_1}=\Ones{D_1}$.  This
corresponds to the fact that if the conditional distribution
within each cloud
is uniform, then taking a random $G_2$-step does nothing.
Hence, $\liftB\ALPHA = \liftB(\ALPHApar+\ALPHAperp) = 
\ALPHApar+\liftB\ALPHAperp$.
Substituting this into (\ref{symmetrize:eqn}), we have
\begin{equation} \label{reflect:eqn}
\langle M\ALPHA,\ALPHA\rangle =
\langle \liftA(\ALPHApar+\liftB\ALPHAperp),\ALPHApar+\liftB\ALPHAperp\rangle.
\end{equation}
Expanding and using the fact that $\liftA$ is 
length-preserving (because it is a permutation matrix),
we have
\begin{eqnarray}
|\langle M\ALPHA,\ALPHA\rangle| 
|\langle \liftA\ALPHApar,\ALPHApar\rangle|
+ 2\|\ALPHApar\|\cdot \|\liftB\ALPHAperp\|
+ \|\liftB\ALPHAperp\|^2.
\label{AB:eqn}
\end{eqnarray}

Now we apply the expansion properties of $G_1$ and $G_2$ to bound
each of these terms.
First, we bound
$\|\liftB\ALPHAperp\|$, which corresponds to 
the intuition that
when the conditional distributions within the clouds are far from
uniform,
they become more uniform when we take a random $G_2$-step.

\numbereddemo{{C}laim} \label{Bbound:clm}
$\|\liftB\ALPHAperp\| \leq \lambda_2 \cdot \|\ALPHAperp\|$.
\enddemo

\demo{Proof of claim}
\begin{eqnarray*}
\liftB\ALPHAperp &=&
\liftB\left(\sum_v e_v\otimes \ALPHA^\perp_v\right)\\
&=& \sum_v e_v\otimes B\ALPHA^\perp_v.
\end{eqnarray*}
By the expansion of $G_2$,
$\|B\ALPHA^\perp_v\|\leq \lambda_2\cdot \|\ALPHA^\perp_v\|$ for all $v$.
Hence, $\|\liftB\ALPHAperp\|\leq \lambda_2\cdot \|\ALPHAperp\|$.
\enddemo

Next, we bound $|\langle \liftA\ALPHApar,\ALPHApar\rangle|$, which
corresponds to the intuition that
when the conditional distribution within each cloud is uniform, the
jump between the clouds makes the marginal distribution
on clouds themselves more uniform.
\numbereddemo{{C}laim} \label{Abound:clm}
$|\langle \liftA\ALPHApar,\ALPHApar\rangle|
\leq \lambda_1\cdot\langle \ALPHApar,\ALPHApar\rangle.$
\enddemo

\demo{Proof of claim}
To prove this, we must first relate $\liftA$ to $A$.
Recall that   when $k$ is uniformly distributed, $\Rot{G_1}(v,k)$
gives a pair $(w,\ell)$ where $w$ is a uniformly selected neighbor
of $v$.  Similarly, if
$\Unit{v}\in\R^{N_1}$ is the $v^{\rm th}$
standard basis vector, then $A\Unit{v}$ gives the uniform distribution
over the neighbors of $v$.  This similarity is captured by the
formula
$\Clouds\liftA(e_v\otimes \Ones{D_1}/D_1)=Ae_v$
for all~$v$.  (Tensoring $e_v$ with $\Ones{D_1}/D_1$ corresponds to
taking the uniform distribution over $k$ and applying
$\Clouds$ corresponds to discarding $\ell$ and looking just at $w$.)
Because the $\Unit{v}$'s form a basis,
this formula extends to all vectors
$\BETA\in\R^{N_1}$:
$\Clouds\liftA(\BETA\otimes \Ones{D_1}/D_1)=A\BETA$.
Applying this formula to
$\ALPHApar=\Clouds\ALPHA\otimes \Ones{D_1}/D_1$, we have
$\Clouds\liftA(\ALPHApar)=A\Clouds\ALPHA$.
Thus,
\begin{eqnarray*}
\langle \liftA\ALPHApar,\ALPHApar\rangle
&=& \langle \liftA\ALPHApar,\Clouds\ALPHA\otimes
\Ones{D_1}\rangle/D_1\\
&=& \langle \Clouds\liftA\ALPHApar,\Clouds\ALPHA\rangle/D_1\\
&=& \langle A\Clouds\ALPHA,\Clouds\ALPHA\rangle/D_1.
\end{eqnarray*}
Recalling that $\Clouds\ALPHA$ is orthogonal to $\Ones{N_1}$, 
we may apply the expansion of $G_1$ to obtain:
\begin{eqnarray*}
|\langle \liftA\ALPHApar,\ALPHApar\rangle|
&\leq&
\lambda_1\cdot \langle \Clouds\ALPHA,\Clouds\ALPHA\rangle/D_1\\
&=&
\lambda_1\cdot \langle \Clouds\ALPHA\otimes \Ones{D_1},
\Clouds\ALPHA\otimes \Ones{D_1}\rangle/D_1^2\\
&=& \lambda_1\cdot \langle \ALPHApar,\ALPHApar\rangle.\\
\noalign{\vskip-36pt}
\end{eqnarray*}
\enddemo
\vglue12pt

Substituting the bounds of Claim~\ref{Bbound:clm} and \ref{Abound:clm}
into (\ref{AB:eqn}), we have:
\begin{equation}\label{cauchyschwartz:eqn}
|\langle M\ALPHA,\ALPHA\rangle| \leq
\lambda_1\cdot\|\ALPHApar\|^2
+ 2\lambda_2\cdot \|\ALPHApar\|\cdot \|\ALPHAperp\|
+ \lambda_2^2\cdot\|\ALPHAperp\|^2.
\end{equation}

If we let $p=\|\ALPHApar\|/\|\ALPHA\|$ and
$q=\|\ALPHAperp\|/\|\ALPHA\|$,
then $p^2+q^2=1$, and the above expression can be rewritten as:
$$
\frac{|\langle M\ALPHA,\ALPHA\rangle|}{\langle \ALPHA,\ALPHA\rangle}
\leq
\lambda_1\cdot p^2
+ 2\lambda_2\cdot pq
+ \lambda_2^2\cdot q^2
\leq 
\lambda_1
+ \lambda_2
+ \lambda_2^2.
$$

This shows that we can take 
$f(\lambda_1,\lambda_2)\leq \lambda_1+\lambda_2+\lambda_2^2$.  It remains to show that we can set 
$f(\lambda_1,\lambda_2)<1$ as long as $\lambda_1,\lambda_2<1$.  
We consider two cases, depending on the length of $\|\ALPHAperp\|$.
First, suppose that 
$\|\ALPHAperp\|\leq \frac{1-\lambda_1}{3\lambda_2}\cdot \|\ALPHA\|.$
Then, from (\ref{cauchyschwartz:eqn}), we have
\begin{eqnarray*}
|\langle M\ALPHA,\ALPHA\rangle| &\leq&
\lambda_1\cdot\|\ALPHA\|^2+ 2\lambda_2\cdot \left(\frac{1-\lambda_1}{3\lambda_2}\right) \|\ALPHA\|^2
\\
&&+\  \lambda_2^2\cdot \left(\frac{1-\lambda_1}{3\lambda_2}\right)^2 \|\ALPHA\|^2
< \left( 1-\frac{1-\lambda_1}{9}\right) \cdot \|\ALPHA\|^2.
\end{eqnarray*}

Now suppose that $\|\ALPHAperp\| > \frac{1-\lambda_1}{3\lambda_2} \cdot \|\ALPHA\|$.
Notice that  $\liftB\ALPHAperp$ is orthogonal to $\ALPHApar$:
$\langle \liftB\ALPHAperp,\ALPHApar\rangle 
= \langle \ALPHAperp,\liftB\ALPHApar\rangle = 
\langle \ALPHAperp,\ALPHApar\rangle = 0.$
Using this, we can bound (\ref{reflect:eqn}) as follows:
\begin{eqnarray*} 
|\langle M\ALPHA,\ALPHA\rangle| &\hskip-6pt =\hskip-6pt&
|\langle \liftA(\ALPHApar+\liftB\ALPHAperp),\ALPHApar+\liftB\ALPHAperp\rangle|
\leq \|\ALPHApar+\liftB\ALPHAperp\|^2
= \|\ALPHApar\|^2+\|\liftB\ALPHAperp\|^2\\
&\hskip-6pt\leq\hskip-6pt& \|\ALPHA\|^2 - \|\ALPHAperp\|^2 + \lambda_2^2\cdot \|\ALPHAperp\|^2
\leq \|\ALPHA\|^2 - 
 (1-\lambda_2^2)\cdot\left(\frac{1-\lambda_1}{3\lambda_2}\right)^2\cdot 
\|\ALPHA\|^2.
\end{eqnarray*}
Thus, we can take
$$f(\lambda_1,\lambda_2)\leq 1-\min\left\{\frac{1-\lambda_1}{9},
\frac{(1-\lambda_1)^2\cdot(1-\lambda_2^2)}{9\lambda_2^2}\right\} < 1.$$


\demo{{\rm 4.2.} Improved analysis of the eigenvalue}
\label{improve-evalue:sec}
In this subsection we state and prove an improved upper bound on the
second largest eigenvalue produced by the zig-zag product. \enddemo

\proclaimtitle{Thm.~\ref{zigzag:thm}, improved}
\proclaim{Theorem} \label{improve-evalue:thm}
If $G_1$ is an $(N_1,D_1,\lambda_1)$\/{\rm -}\/graph and $G_2$ is
a $(D_1,D_2,\lambda_2)$\/{\rm -}\/graph{\rm ,} then
$G_1\zigzag G_2$ is a
$(N_1\cdot D_1,D_2^2,
f(\lambda_1,\lambda_2))$\/{\rm -}\/graph{\rm ,}
where
$$f(\lambda_1,\lambda_2)=\frac{1}{2}(1-\lambda_2^2)\lambda_1+
\frac{1}{2}\sqrt{(1-\lambda_2^2)^2\lambda_1^2+4\lambda_2^2}.$$
\endproclaim

Although the function $f(\lambda_1,\lambda_2)$ looks ugly, it can be verified
that it has the following nice properties:
\begin{itemize}
\item[1.] $f(\lambda,0)=f(0,\lambda)=\lambda$ and $f(\lambda,1)=f(1,\lambda)=1$ for all $\lambda\in [0,1]$.
\item[2.] $f(\lambda_1,\lambda_2)$ is a strictly increasing function of both $\lambda_1$ and $\lambda_2$ (except when
one of them is 1).
\item[3.] If $\lambda_1<1$ and $\lambda_2<1$, then $f(\lambda_1,\lambda_2)<1$.
\item[4.] $f(\lambda_1,\lambda_2)\leq \lambda_1+\lambda_2$ for all $\lambda_1, \lambda_2 \in [0,1]$.
\end{itemize}

\demo{{P}roof}
The proof proceeds along the same lines as
the proof of Theorem~\ref{zigzag:thm}, except that we will use a
geometric argument
to directly bound (\ref{reflect:eqn}) rather than first passing
to (\ref{AB:eqn}). That is, we
must bound (using the same notation as in that proof)
$$
\frac{\langle M\ALPHA,\ALPHA\rangle}{\langle \ALPHA,\ALPHA\rangle} =
\frac{\langle \liftA(\ALPHApar+\liftB\ALPHAperp),\ALPHApar+\liftB\ALPHAperp\rangle}{\|\ALPHApar+\ALPHAperp\|^2}.
$$
The key observation is:
\numbereddemo{{C}laim} \label{reflect:clm}
$\liftA$ is a reflection through a linear subspace
$S$ of $\R^{N_1D_1}$.  Hence, for any vector $v$,
$\langle \liftA v,v\rangle = ({\rm cos} 2\theta)\cdot \|v\|^2,$ where
$\theta$ is the angle between $v$ and~$S$.
\enddemo

{\it Proof of claim}.
By the symmetry of $\liftA$, we can decompose $\R^{N_1D_1}$ into 
the sum of orthogonal eigenspaces of $\liftA$.
Since
$\liftA^2=I_{N_1D_1}$, the only eigenvalues of $\liftA$ are $\pm 1$.
Take $S$ to be the $1$-eigenspace of $\liftA$.
\hfill\qed
\vglue12pt

Thus,
the expression we want to bound is
$$\frac{|\langle M\ALPHA,\ALPHA\rangle |}{\langle \ALPHA,\ALPHA\rangle} = 
|{\rm cos} 2\theta|\cdot 
\frac{\|\ALPHApar+\liftB\ALPHAperp\|^2}{\|\ALPHApar+\ALPHAperp\|^2}
=
|{\rm cos} 2\theta|\cdot 
\frac{{\rm cos}^2 \phi}{{\rm cos}^2 \phi'},$$
where $\theta$ is the angle between $\ALPHApar+\liftB\ALPHAperp$
and $S$, $\phi\in[0,\pi/2]$ is the angle between $\ALPHApar$ and 
$\ALPHApar+\ALPHAperp$, and $\phi'\in[0,\pi/2]$ is the angle between
$\ALPHApar$ and $\ALPHApar+\liftB\ALPHAperp$.
If we also let $\psi$ be the angle between $\ALPHApar$ and $S$,
then we clearly have $\theta \in [\psi-\phi',\psi+\phi']$.

Now we translate Claims~\ref{Bbound:clm} and \ref{Abound:clm} into
this geometric language.  Claim~\ref{Bbound:clm} constrains
the relationship between $\phi'$ and $\phi$ by 
$$\frac{\tan \phi'}{\tan \phi}=\frac{\|\liftB\ALPHAperp\|}{\|\ALPHAperp\|}\leq \lambda_2.$$
Claim~\ref{Abound:clm} says
$|{\rm cos} 2\psi|\leq \lambda_1$.  
For notational convenience, we will denote the exact values of
$(\tan \phi')/(\tan \phi)$ 
and  $|{\rm cos} 2\psi|$  by $\mu_2$ and $\mu_1$, respectively.  We will
work with these values until the end of the proof, at which point
we will upper bound them by $\lambda_2$ and $\lambda_1$.

To summarize, we want to maximize
\begin{equation} \label{cos:eqn}
|{\rm cos} 2\theta|\cdot 
\frac{{\rm cos}^2 \phi}{{\rm cos}^2 \phi'} 
\end{equation}
over the variables $\theta$, $\phi$, $\phi'$, and $\psi$, subject to the
following constraints:
\begin{itemize}
\item[1.] $\phi,\phi',\psi\in [0,\pi/2]$.
\item[2.] $\theta\in [\psi-\phi',\psi+\phi']$.\footnote{We
do not require $\theta\in [0,\pi/2]$ so that we do not have to worry
about ``wraparound''  in the interval $[\psi-\phi',\psi+\phi']$.  
Adding a multiple of $\pi/2$ to $\theta$ does not change the value
of (\ref{cos:eqn}).}
\item[3.] $\tan \phi'/\tan \phi= \mu_2$.
\item[4.] $|{\rm cos} 2\psi| = \mu_1$.
\end{itemize}
There are two cases, depending on whether 
$|{\rm cos} 2x|$ ever achieves the value 1 in the interval
$[\psi-\phi',\psi+\phi']$.

\demo{Case I} $\phi' \leq \min\{\psi,\pi/2-\psi\}$. 
Then
\begin{eqnarray*}
|{\rm cos} 2\theta| &=& \max\{|{\rm cos} 2(\psi+\phi')|,|{\rm cos} 2(\psi-\phi')|\}\\
        &=& |{\rm cos} 2\psi\cdot{\rm cos} 2\phi'|+|\sin 2\psi\cdot\sin 2\phi'|.
\end{eqnarray*}
After some trigonometric manipulations, we have
$$|{\rm cos} 2\theta|\cdot \frac{{\rm cos}^2 \phi}{{\rm cos}^2 \phi'}  
\!= \!\frac{1}{2}\left|(1-\mu_2^2){\rm cos} 2\psi+(1+\mu_2^2){\rm cos} 
2\psi{\rm cos} 2\phi\right|+\frac{1}{2}|2\mu_2\sin 2\psi\sin 2\phi|.$$
The choice of $\phi$ which maximizes this is to have 
$({\rm cos} 2\phi,\sin 2\phi)$ be a unit vector in the direction
of $(\pm (1+\mu_2^2){\rm cos} 2\psi,2\mu_2\sin 2\psi)$; thus
\begin{eqnarray*}
|{\rm cos} 2\theta|\cdot \frac{{\rm cos}^2 \phi}{{\rm cos}^2 \phi'}  
&\leq & \frac{1}{2}(1-\mu_2^2)|{\rm cos} 2\psi|+
\frac{1}{2}\sqrt{(1+\mu_2^2)^2{\rm cos}^2 2\psi+4\mu_2^2\sin^2 2\psi}\\
&=& \frac{1}{2}(1-\mu_2^2)\mu_1+
\frac{1}{2}\sqrt{(1+\mu_2^2)^2\mu_1^2+4\mu_2^2(1-\mu_1^2)}.\\
\noalign{\vskip-30pt}
\end{eqnarray*}
\enddemo

\demo{Case II} $\phi' > \min\{\psi,\pi/2-\psi\}$. 
In this case, we cannot obtain any nontrivial bound
on $|{\rm cos} 2\theta|$, so, after some
trigonometric manipulations, the problem is reduced to bounding:
\begin{equation} \label{case2:eqn}
|{\rm cos} 2\theta|\cdot \frac{{\rm cos}^2 \phi}{{\rm cos}^2 \phi'}  
\leq \frac{{\rm cos}^2 \phi}{{\rm cos}^2 \phi'}  = \mu_2^2+(1-\mu_2^2){\rm cos}^2 \phi.
\end{equation}
The condition $\phi' > \min\{\psi,\pi/2-\psi\}$
implies that ${\rm cos} 2\phi' < |{\rm cos} 2\psi| = \mu_1$.
After some trigonometric manipulations, we have
$${\rm cos} 2\phi' = \frac{(1+\mu_2^2){\rm cos}^2 \phi-\mu_2^2}{(1-\mu_2^2){\rm cos}^2 \phi+\mu_2^2},$$
and the condition ${\rm cos} 2\phi' < \mu_1$ is equivalent to
$${\rm cos}^2 \phi < \frac{\mu_2^2(1+\mu_1)}{(1-\mu_1)+\mu_2^2(1+\mu_1)}.$$
Substituting this into (\ref{case2:eqn})
and simplifying, we
conclude that
$$|{\rm cos} 2\theta|\cdot \frac{{\rm cos}^2 \phi}{{\rm cos}^2 \phi'}<
\frac{2\mu_2^2}{1-\mu_1+\mu_2^2(1+\mu_1)}.$$

It can be verified that the bound obtained in Case I is an
increasing function of $\mu_1$ and $\mu_2$ and is always
greater than or equal to the bound in Case~II.  Therefore,
replacing $\mu_1$ and $\mu_2$ by $\lambda_1$ and $\lambda_2$
in the Case I bound proves the theorem.
\hfill\qed
\enddemo

\vglue-16pt

\section{The base graph} \label{basegraph:sec}
\vglue-6pt

Our construction of an infinite family of expanders
in Section~\ref{recursion:sec} requires starting with
a $(D^8,D,\lambda)$-graph $H$ 
(for
a sufficiently small $\lambda$, say $\leq 1/5$).
Since $D$ is a ``constant,'' such a graph can be found by
exhaustive search (given that one exists, which can be
proven by (nontrivial) probabilistic 
arguments~\cite{Alo86a}, \cite{BS87}, \cite{FKS89}, \cite{Fri91}).  
However, for these parameters,
there are simple explicit constructions known.
We describe two of them below.
The first is simpler and more intuitive, but
the second yields better parameters.

\demo{{\rm 5.1.} The affine plane}
The first construction is based on the ``projective plane''
construction of Alon~\cite{Alo86b}, but we instead use the affine plane
in order to make 
$N$ exactly $D^2$ and then use
the zig-zag product to obtain a graph with $N=D^8$.
For a prime power $q=p^t$, let $\F_q$ be the finite field
of size $q$; an explicit representation of such a field
can be found deterministically in time $\poly(p,t)$~\cite{Sho90}.  
We define a graph
$\AP{q}$ with vertex set $\F_q^2$, and edge set
$\{((a,b),(c,d)) : ac=b+d\}$.  
That is, we connect 
the vertex $(a,b)$ to all points on the line $L_{a,b}=\{(x,y) : y=ax-b\}$.
(Note that we have chosen the sign of $b$ to make the graph
undirected.)\enddemo

\proclaim{Lemma}
$\AP{q}$ is an $(q^2,q,1/\sqrt{q})${\rm -}graph{\rm .}
Moreover{\rm ,} a rotation map for $\AP{q}$ can be computed in
time $\poly(\log q)$ given a representation of the field~$\F_q${\rm .}
\endproclaim

{\it Proof}.
The expansion of $\AP{q}$ will follow from the fact the square
of $\AP{q}$ is almost the complete graph, which in turn
is based on the fact that almost all pairs of lines in the
plane $\F_q^2$ intersect.  Let $M$ be the
$q^2\times q^2$ normalized adjacency matrix of $\AP{q}$; 
we will now calculate the
entries of $M^2$.  The entry of $M^2$ in row $(a,b)$ and
column $(a',b')$ is
exactly the number of common neighbors of $(a,b)$ and $(a',b')$
in $\AP{q}$ divided by $q^2$, i.e., $|L_{a,b}\cap L_{a',b'}|/q^2$.
If $a\neq a'$, then $L_{a,b}$ and $L_{a',b'}$ intersect in exactly
one point.  If $a=a'$ and $b\neq b'$, then their intersection is
empty, and if $a=a'$ and $b=b'$, then their intersection is of size
$q$.  Thus, if we let $I_q$ denote the $q\times q$ identity matrix
and $J_q$ the $q\times q$ all-1's matrix, we have
$$M^2 = \frac{1}{q^2} \left( \begin{array}{cccc}
qI_q & J_q & \cdots & J_q \\ 
J_q & qI_q &  & J_q \\
\vdots &  & \ddots & J_q\\
J_q & J_q & \cdots & qI_q
\end{array} \right) = \frac{I_q \otimes qI_q + (J_q-I_q)\otimes J_q}{q^2}.$$
Now we can calculate the eigenvalues explicitly.
$J_q$ has eigenvalues $q$ (multiplicity 1) and $0$ (multiplicity $q-1$).  
So $(J_q-I_q)\otimes J_q$ has eigenvalues
$(q-1)\cdot q$, $-1\cdot q$, and $0$.  Adding $I_q\otimes qI_q$
increases all these eigenvalues by $q$, and then we divide by $q^2$.
Hence the eigenvalues of $M^2$ are 1 (multiplicity 1), 0 
(multiplicity $q-1$), and $1/q$ (multiplicity $(q-1)\cdot q$).
Therefore, the
second largest eigenvalue of $M$ has absolute value $1/\sqrt{q}$.

A rotation map for $\AP{q}$ is given by 
$$\Rot{q}((a,b),t)=\cases{
((t/a,t-b),t) & if $a\neq 0$ and $t\neq 0$,\cr
((t,-b),a) & if $a=0$ or $t=0$,}$$
where $a,b,t\in \F_q$.
\hfill\qed\vglue9pt

Now, define the following graphs inductively:
\begin{eqnarray*}
\AP{q}^1 &=& \AP{q}\otimes \AP{q}\\
\AP{q}^{i+1} &=& \AP{q}^i \zigzag \AP{q}.
\end{eqnarray*}
From Proposition~\ref{tensor:prop} and Theorem~\ref{zigzag:thm},
we immediately deduce:

\proclaim{Proposition} \label{linepoint:prop}
$\AP{q}^i$ is a $(q^{2(i+1)},q^2,O(i/\sqrt{q}))$\/{\rm -}\/graph.\footnote{The
hidden constant in $O(i/\sqrt{q})$ can be reduced to $1$
using the improved analysis of the zig\/{\rm -}\/zag product in
Theorem~\ref{improve-evalue:thm}.}
Moreover{\rm ,} a rotation map for $\AP{q}^i$ can be computed in
time $\poly(i,\log q)$ given a representation of $\F_q${\rm .}
\endproclaim

Taking $i=7$ and a sufficiently large $q$ gives a graph
suitable for the expander construction in Section~\ref{recursion:sec}.

\vglue6pt 5.2. {\it Low\/{\rm -}\/degree polynomials}.  
The graphs we describe here are derived
from constructions of Alon and Roichman~\cite{AR94}, which are
Cayley graphs derived from
the generator matrix of an error-correcting code.
In order to give a self-contained presentation, we 
specialize the construction to a Reed-Solomon code concatenated with a
Hadamard code (as used in, e.g. \cite{AGHP92}).

For a prime power $q$ and $d\in \N$, we define
a graph $\LD_{q,d}$ on vertex set $\F_q^{d+1}$ with degree $q^2$.
For a vertex $a\in\F_q^{d+1}$ and $x,y\in \F_q$, the
the $(x,y)^{\rm th}$ neighbor of $a$ is
$a+(y,yx,yx^2,\ldots,yx^d)$. 

\proclaim{Proposition} \label{polygraph:prop}
$\LD_{q,d}$ is a $(q^{d+1},q^2,d/q)$\/{\rm -}\/graph{\rm .}  Moreover{\rm ,} a
rotation map for $\LD_{q,d}$ can be computed in time 
$\poly(\log q,d)$ given a representation of~$\F_q${\rm .}
\endproclaim

As above, taking $d=7$ and sufficiently large $q$
gives a graph suitable for our expander construction.  These
graphs are better than those of Proposition~\ref{linepoint:prop} because the
the eigenvalue-degree relationship is the optimal $\lambda=O(1/\sqrt{D})$
(as $q$ grows), which implies an eigenvalue of $O(1/D^{1/4})$ for the
family constructed in Theorem~\ref{expfam:thm}.

\demo{Proof}
To simplify notation, let $\F=\F_q$.
Let $M$ be the $q^{d+1}\times q^{d+1}$ 
normalized adjacency matrix of $\LD_{q,d}$.
We view vectors in $\C^{q^{d+1}}$ as functions
$ f : \F^{d+1}\to \C$.
We will now explicitly describe the eigenvectors
of $M$.  Let $p$ be the characteristic of $\F$, let
$\zeta=e^{2\pi i/p}$ be a primitive $p^{\rm th}$ root of unity,
and let $L :\F\to \F_p$ be any surjective
$\F_p$-linear map. (For simplicity, 
one can think of the special case that $p=q$ and
$L$ is the identity map.)

For every sequence $a=(a_0,\ldots,a_d)\in \F^{d+1}$, define
the function $\chi_a : \F^{d+1} \to \C$
by $\chi_a(b) = \zeta^{L(\sum a_ib_i)}$.
Clearly,
$\chi_a(b+c)=\chi_a(b)\chi_a(c)$
for any $b,c\in\F^{d+1}$.
Moreover, it can be verified that the $\{\chi_a\}$ are
orthogonal under the standard
inner product $\langle f,g\rangle = \sum_b f(b)g(b)^*$,
and thus form a basis for $\C^{q^{d+1}}$.
Hence, if we show that each $\chi_a$ is an eigenvector
of $M$, then they are all the eigenvectors of $M$.
This can be done by direct calculation:

\begin{eqnarray*}
(M\chi_a)(b) &=& \frac{1}{q^2}\sum_{c\in \F^{d+1}} M_{bc}\cdot\chi_a(c)\\
&=& \frac{1}{q^2} \sum_{x,y\in \F} \chi_a(b+(y,yx,\ldots,yx^d))\\
&=& \left(\frac{\sum_{x,y\in\F} \chi_a(y,yx,\ldots,yx^d)}{q^2}\right)\cdot
\chi_a(b)\\
&\eqdef& \lambda_a \cdot\chi_a(b).
\end{eqnarray*}

Thus, $\chi_a$ is an eigenvector of $M$ with eigenvalue $\lambda_a$
and all eigenvectors of $M$ are of this form.  
So we simply need
to show that $|\lambda_a|\leq d/q$ for all but one $a\in\F^{d+1}$.
To do this, note that 
$$\lambda_a=\frac{1}{q^2}\sum_{x,y\in\F}
\chi_a((y,yx,\ldots,yx^d))=\frac{1}{q^2}\sum_{x,y\in\F}
\zeta^{L(y\cdot p_a(x))},$$
where $p_a(x)$ is the polynomial $a_0+a_1x+\cdots+a_dx^d$.
When $x$ is a root of $p_a$, then $\zeta^{L(yp_a(x))}=1$ for all $y$;  hence $x$ contributes $q/q^2=1/q$ to $\lambda_a$.
When $x$ is not a root of $p_a(x)$, $yp_a(x)$ takes on all values
in $\F$ as $y$ varies, and hence $\zeta^{L(yp_a(x))}$ varies
uniformly over all $p^{\rm th}$ roots of unity.  Since the sum of
all $p^{\rm th}$ roots of unity is 0, these $x$'s contribute nothing
to $\lambda_a$.
When $a\neq 0$, $p_a$ has at most $d$ roots,
so $|\lambda_a|\leq d/q$.
\enddemo

\section{Variants on the zig-zag theme} \label{zigzagmod:sec}

The two subsections of this section contain two variants of the basic
zig-zag product. The first is aimed at 
improving the relation between the degree and the eigenvalue bound. 
The second is aimed at simplifying
the product, at the cost of deteriorating this relationship.
 
\demo{{\rm 6.1.} A {\rm ``}\/derandomized\/{\rm ''} zig\/{\rm -}\/zag product} 
In this section we provide a variant of our original zig-zag product,
which achieves a better relationship between the degree and the
expansion of the resulting graph. The term ``derandomized'' will
become clearer when we define it.

Recall that the optimal second-largest eigenvalue 
for an infinite family
of $D$-regular graphs is $\Theta(1/D^{1/2})$, and families of
graphs meeting this bound (with the right constant)
are referred to as Ramanujan. A basic question is how close can we come
to this optimal bound using our techniques.
Starting with a constant-size Ramanujan graph
(or the graphs of \S 5.2),
our basic construction
of Theorem~\ref{expfam:thm} achieves a second-largest eigenvalue
of $O(1/D^{1/4})$ for the family of expanders generated..  

Here, we define a variant of the zig-zag product, which makes more
efficient use of the expansion of the small graph.
Using the new product
in our iterative construction
(of \S \ref{recursion:sec})
with an initial constant-size Ramanujan graph
or even the graphs of Proposition~\ref{polygraph:prop}, 
we obtain a second-largest
eigenvalue of $O(1/D^{1/3})$ for the family of expanders generated.
It is an interesting open problem to construct
families of graphs achieving the optimal
eigenvalue $O(1/D^{1/2})$ using a similar graph product.

We now turn to the formal definition of the new zig-zag product. It
will have two ``zig'' moves and two ``zag'' moves, but they will not
be independent. The second ``zig'' and the first ``zag'' will use the
same random bits!

\numbereddemo{Definition}\label{zigzagmod:dfn}
Let
$G_1$ be a $D_1$-regular graph on $[N_1]$ with rotation
map $\Rot{G_1}$ and let $G_2$ be a $D_2$-regular graph on $[D_1]$
with rotation map $\Rot{G_2}$.
Suppose that 
for every $i\in [D_2]$, $\Rot{G_2}(\cdot,i)$ induces
a permutation on $[D_1]$.\footnote{By this we mean that 
the function $f_i(x)=\mbox{``the first component of $\Rot{G_2}(x,i)$''}
= \mbox{``the $i^{\rm th}$ neighbor of $x$''}$ is a permutation
for every $i$.}
Then the 
{\it modified zig-zag product}
of $G_1$ and $G_2$ 
is defined to be the $D_2^3$-regular graph 
$G_1\zigzagmod G_2$ on $[N_1]\times [D_1]$
whose rotation map $\Rot{G_1\zigzagmodsub G_2}$ is as follows:
\medbreak {\it $\Rot{G_1\protect\zigzagmodsub G_2}((v,k),(h,i,j))$}:
\begin{itemize}
\item[1.] Let $(k',h')=\Rot{G_2}(k,h)$.
\item[2.] Let $(k'',i')=\Rot{G_2}(k',i)$. 
\item[3.] Let $(w,\ell'')=\Rot{G_1}(v,k'')$.
\item[4.] Find the unique $\ell'\in [D_1]$ such that $(\ell'',i'')=\Rot{G_2}(\ell',i)$ for some $i''$. ($\ell'$ exists by the
assumption on $\Rot{G_2}$.)
\item[5.] Let $(\ell,j')=\Rot{G_2}(\ell',j)$. 
\item[6.] Output $((w,\ell),(j',i,h'))$.
\end{itemize}
\enddemo

Again, in this graph 
product we do {\it two} random steps on the small graph in
both the zig and the zag parts. However, to save random bits (i.e.,
decrease the degree) we use {\it the same} random bits for the second
move of the zig part and the first move of the zag part. Thus the
degree of the new graph is $D_2^3$. However, we will show that the
bound on the eigenvalue will be as if these moves were independent.
This proof will follow the lines of the basic analysis of the original
zig-zag product.

\proclaim{Theorem}\label{zigzagmod:thm}
If $G_1$ is an $(N_1,D_1,\lambda_1)$\/{\rm -}\/graph and $G_2$ is
a $(D_1,D_2,\lambda_2)$\/{\rm -}\/graph{\rm ,} then
$G_1\zigzagmod G_2$ is a
$(N_1\cdot D_1,D_2^3,\lambda_1 + 2\lambda_2^2 )$\/{\rm -}\/graph{\rm .}
Moreover{\rm ,} $\Rot{G_1\zigzagmodsub G_2}$ can be computed in time
$\poly(\log N,\log D_1,D_2)$ with one oracle query to
$\Rot{G_1}$ and $D_2+2$ oracle queries to $\Rot{G_2}${\rm .}
\endproclaim

{\it Proof}.
We use the same notation as in the proof of 
Theorem~\ref{zigzag:thm}.  Like there, we need to bound 
$|\langle M\ALPHA,\ALPHA\rangle|/\langle \ALPHA,\ALPHA\rangle$,
where $M$ is the normalized adjacency matrix of $G_1\zigzagmod G_2$
and $\ALPHA \perp \Ones{N_1D_1}$.  Let $B_i$ be the $D_1\times D_1$
permutation matrix induced by $\Rot{G_2}(\cdot,i)$, and 
let $\liftB_i=I_{N_1}\otimes B_i$.  Then 
$$\liftB=\frac{1}{D_2}\sum_{i=1}^{D_2} \liftB_i.$$  
Note that the normalized adjacency matrix 
corresponding to Steps~2--4 in
the definition of $G_1\zigzagmod G_2$ is given by
$$\Mmod=\frac{1}{D_2} \sum_i \liftB_i\liftA\liftB_i^T,$$
where $\liftB_i^T$ is the transpose (equivalently, inverse) of
$\liftB_i$.
Thus, $M=\liftB\Mmod\liftB$.
The main observation is that not only does $\liftB\ALPHApar=\ALPHApar$
(as we used in the original analysis), but also $\liftB_i^T\ALPHApar=\ALPHApar$
for every $i$ (because $B_i$ is a permutation matrix).
Hence,
$$ \Mmod\ALPHApar = \frac{1}{D_2} \sum_i \liftB_i\liftA\liftB_i^T\ALPHApar
= \frac{1}{D_2} \sum_i \liftB_i\liftA\ALPHApar
= \liftB\liftA\ALPHApar.$$
Applying this (and the symmetry of $\liftB$ and $\Mmod$), we get
\begin{eqnarray*}
\langle M\ALPHA,\ALPHA\rangle &=&
\langle M\ALPHApar,\ALPHApar\rangle +
2\langle M\ALPHApar,\ALPHAperp\rangle+
\langle M\ALPHAperp,\ALPHAperp\rangle\\[5pt]
&=& \langle \liftA\ALPHApar,\ALPHApar\rangle
+ 2\langle \liftA\ALPHApar,\liftB^2\ALPHAperp\rangle
+
\langle \Mmod\liftB\ALPHAperp,\liftB\ALPHAperp\rangle.
\end{eqnarray*}
Being the normalized adjacency matrix of an undirected, regular graph,
$\Mmod$ has no eigenvalues larger than 1 and hence
does not increase the length of any vector.
Using this together with
 Claims~\ref{Bbound:clm} and \ref{Abound:clm},
we have
\begin{eqnarray*}
|\langle M\ALPHA,\ALPHA\rangle| &\leq&
|\langle \liftA\ALPHApar,\ALPHApar\rangle|
+ 2\|\ALPHApar\|\cdot\|\liftB^2\ALPHAperp\|
+ \|\liftB\ALPHAperp\|^2\\[5pt]
&\leq& \lambda_1\cdot\|\ALPHApar\|^2
+ 2\lambda_2^2\cdot \|\ALPHApar\|\cdot \|\ALPHAperp\|
+ \lambda_2^2\cdot\|\ALPHAperp\|^2.
\end{eqnarray*}
As in the the proof of Theorem~\ref{zigzag:thm},
using the fact that $\|\ALPHApar\|^2+\|\ALPHAperp\|^2=\|\ALPHA\|^2$
yields the desired bound.
\hfill\qed

\demo{{\rm 6.2.} The replacement product}
In this section, we describe an extremely simple and intuitive graph
product, which shares similar properties to the zig-zag
product. Namely, when taking the product of two expanders, we get
a larger expander whose degree depends only on that of the smaller
graph. Here simplicity is the important feature, and the expansion
quality is not as good as above.  This product is so
natural that it was used in various contexts before. 
Indeed, Gromov \cite{Gro83} even estimates the second eigenvalue of an
iterated replacement product of the graph of the Boolean hypercube with
smaller copies of itself. (Of course, in this very special case the
outcome is not expanding, since the cube is not.)
Our proof of its expansion will be  a simple reduction to
the expansion properties of the zig-zag product. However, one can
also prove it directly in a manner similar to the proof of
Theorem~\ref{zigzag:thm} (and thereby obtain a stronger bound). 

 Assume (as in the basic zig-zag product) that 
$G_1$ is a $D_1$-regular graph on $[N_1]$ and 
$G_2$ is a $D_2$-regular graph on $[D_1]$.
A natural idea is to place a ``copy'' (or ``cloud'') 
of $G_2$ around each vertex of
$G_1$, maintaining the edges of both. More precisely, every vertex
will be connected to all its original neighbors in its cloud, as well
as to one vertex in the neighboring cloud it defines.
For example, if $G_1$ is
the $n$-dimensional Boolean cube graph, and $G_2$ is the cycle on $n$
vertices, then the resulting graph is the so-called
{\it cube connected cycle}, 
which used to be a popular architecture
for parallel computers.
Note that in this example the small graph had degree 2, and the product
graph had degree~3. In general, the resulting graph would have degree
$D_2+1$. 
In terms of rotation maps, this product is defined as follows.
\numbereddemo{Definition}If $G_1$ is a $D_1$-regular graph on $[N_1]$ with rotation
map $\Rot{G_1}$ and $G_2$ is a $D_2$-regular graph on $[D_1]$
with rotation map $\Rot{G_2}$, then their 
{\it replacement product}
$G_1\replacement G_2$ 
is defined to be the $(D_2+1)$-regular graph on $[N_1]\times [D_1]$ 
whose rotation map $\Rot{G_1\replacementsub G_2}$ is as follows:
\medbreak {\it $\Rot{G_1\protect\replacementsub G_2}((v,k),i)$}:
\vglue4pt
 1.  If $i\leq D_2$, let $(m,j)=\Rot{G_2}(k,i)$ and output
$((v,m),j)$.
\vglue4pt 2.  If $i=D_2+1$, output $(\Rot{G_1}(v,k),i)$.
\enddemo

The expansion properties of the replacement product are given in the
next theorem, relating it to those of the zig-zag product.

\proclaim{Theorem} \label{replacement:thm}
If $G_1$ is an $(N_1,D_1,\lambda_1)$\/{\rm -}\/graph and $G_2$ is
a $(D_1,D_2,\lambda_2)$\/{\rm -}\/graph{\rm ,} then
$G_1\replacement G_2$ is a
$(N_1\cdot D_1,D_2+1,g(\lambda_1,\lambda_2,D_2))$\/{\rm -}\/graph{\rm ,}
where (using the function $f$ from Theorems~{\rm \ref{zigzag:thm}}
or {\rm \ref{improve-evalue:thm})}
$$g(\lambda_1,\lambda_2,D_2) \leq \left(p+(1-p)f(\lambda_1,\lambda_2)\right)^{1/3},$$
and $p=D_2^2/(D_2+1)^3${\rm .} In particular{\rm ,} 
$g(\lambda_1,\lambda_2,D_2)<1$ when $\lambda_1,\lambda_2<1${\rm .}
Moreover{\rm ,} $\Rot{G_1\replacementsub G_2}$ can be computed in time
$\poly(\log N,\log D_1,\log D_2)$ with one oracle query to
$\Rot{G_1}$ or $\Rot{G_2}${\rm .}
\endproclaim

{\it Proof}.
The idea of the proof is that the graph of the zig-zag product is
a regular subgraph of {\it the cube} of the graph of the replacement
product. 
Let $M$ denote the normalized adjacency matrix of $G_1
\replacement G_2$. As in the proof of Theorem~\ref{zigzag:thm}, we let
$A,B$ respectively denote the normalized adjacency matrices of
$G_1,G_2$, and define their ``liftings'' $\liftA , \liftB $ in the
same way. By inspection, we have $M= ( \liftA + D_2\liftB )/(D_2+1)$. 
The key observation is that
$$M^3 = \frac{(\liftA+D_2\liftB)^3}{(D_2+1)^3} = 
p\liftB\liftA\liftB
+(1-p)C,$$
where  $\liftB \liftA \liftB $ is the normalized adjacency matrix of
$G_1 \zigzag\, G_2$, $C$ is the normalized adjacency
matrix of an undirected, regular graph (and in particular
does not increase the length of any vector),
and $p=D_2^2/(D_2+1)^3$.
As eigenvalues of powers of matrices are the respective powers of the
original eigenvalues (see Proposition~\ref{power}), we have
\vglue4pt
\hfill $g(\lambda_1,\lambda_2) \leq (p+ (1-p)f(\lambda_1,\lambda_2))^{1/3}.$ 
\hfill\qed
\pagebreak

Thus, for ``constant'' degrees $D_2$ the replacement product indeed
transforms two expanders into a larger one. 
As in Corollary~\ref{deg4:cor}, we can use this to get 
degree 3 expanders.

\proclaim{{C}orollary}\label{deg3:cor}
For every $\lambda<1$ and every odd $D${\rm ,} there exists a $\lambda'<1$ such that
if $G$ is an $(N,D,\lambda)$\/{\rm -}\/graph and $C$ is the cycle on $D$ vertices{\rm ,}
then $G\replacement C$ is a $(ND,3,\lambda')$\/{\rm -}\/graph{\rm .} 
\endproclaim 

To make the expansion properties in Theorem~\ref{replacement:thm} 
independent of how large $D_2$ is, we now
slightly modify the replacement product to have $D_2$
copies of each edge which goes between clouds. This makes the degree
of every vertex $2D_2$, of which $D_2$ stay within the same cloud, and
the other $D_2$ all connect to the same vertex in a neighbor cloud.
This ``balancing'' make the random walk give the same weight to edges
defined by $G_1$ and $G_2$.

\numbereddemo{Definition} If $G_1$ is a $D_1$-regular graph on $[N_1]$ with rotation
map $\Rot{G_1}$ and $G_2$ is a $D_2$-regular graph on $[D_1]$
with rotation map $\Rot{G_2}$, then their {\it balanced
replacement product}
$G_1\balanced G_2$ 
is defined to be the $2D_2$-regular graph on $[N_1]\times [D_1]$ 
whose rotation map $\Rot{G_1\balancedsub G_2}$ is as follows:
\medbreak {\it $\Rot{G_1\protect\balancedsub G_2}((v,k),i)$}:
\begin{itemize}
\item[1.] If $i\leq D_2$, let $(m,j)=\Rot{G_2}(k,i)$ and output
$((v,m),j)$.
\item[2.] If $i>D_2$, output $(\Rot{G_1}(v,k),i)$.
\end{itemize}
\enddemo

\proclaim{Theorem} \label{balanced:thm}
If $G_1$ is an $(N_1,D_1,\lambda_1)$\/{\rm -}\/graph and $G_2$ is
a $(D_1,D_2,\lambda_2)$\/{\rm -}\/graph{\rm ,} then
$G_1\balanced G_2$ is a
$(N_1\cdot D_1,2D_2,h(\lambda_1,\lambda_2))$\/{\rm -}\/graph{\rm ,}
where {\rm (}\/using the function $f$ from Theorems {\rm \ref{zigzag:thm}}
or {\rm \ref{improve-evalue:thm})}
$$h(\lambda_1,\lambda_2) \leq \left(\frac78+ \frac18\cdot f(\lambda_1,\lambda_2)\right)^{1/3}.$$
In particular{\rm ,} 
$h(\lambda_1,\lambda_2)<1$ when $\lambda_1,\lambda_2<1${\rm .}
Moreover{\rm ,} $\Rot{G_1\balancedsub G_2}$ can be computed in time
$\poly(\log N,\log D_1,\log D_2)$ with one oracle query to
$\Rot{G_1}$ and one oracle query to $\Rot{G_2}$.
\endproclaim

\demo{Proof}
The proof is the same as that of Theorem~\ref{replacement:thm},
noting instead that $M=(\liftA+\liftB)/2$.
\enddemo

As a final note, we observe the weakness of the replacement products
relative to the zig-zag product. Informally, in zig-zag the expansion
quality of the product improves with those of its component, while in
the replacement it does not. More formally,
while the function
$f(\lambda_1,\lambda_2)$ tends to zero when $\lambda_1$ and $\lambda_2$ do, the
functions $g(\lambda_1,\lambda_2,D_2)$ and $h(\lambda_1,\lambda_2)$ do not.

\demo{Acknowledgments}
We are grateful to David Zuckerman for illuminating discussions
and a number of useful suggestions early in the stages of this
work.
We thank the organizers of the
DIMACS Workshop on Pseudorandomness and Explicit Combinatorial
Constructions in October 1999, where we began this research.
We are grateful to Peter Winkler for suggesting the name
``zig-zag product.''
We also thank Neil Agarwol, Noga Alon,
Oded Goldreich, Peter Sarnak, Ronen Shaltiel, 
Dan Spielman, and the anonymous referee for helpful comments
and pointers.
\enddemo

\AuthorRefNames [RVWOOO]

\end{document}